\documentclass[12pt,reqno]{amsart}
\usepackage{amsthm,amsfonts,amssymb,amscd}
\usepackage[all]{xy}

\newtheorem{theorem}[subsection]{Theorem}

\newtheorem{proposition}[subsection]{Proposition}

\newtheorem{lemma}[subsection]{Lemma}
\newtheorem{corollary}[subsection]{Corollary}

\theoremstyle{definition}

\newtheorem{proposition-definition}[subsection]{Proposition-Definition}

\theoremstyle{remark}
\newtheorem{remark}[subsection]{Remark}

\newcommand{\dual}{{\scriptscriptstyle \vee}}

\newcommand{\codim}{\operatorname{codim}\nolimits}

\newcommand{\rk}{\operatorname{rk}\nolimits}

\newcommand{\Homg}{\operatorname{Hom}\nolimits}
\newcommand{\id}{\operatorname{id}\nolimits}
\newcommand{\Ext}{\operatorname{Ext}\nolimits}
\newcommand{\Alb}{\operatorname{Alb}\nolimits}
\newcommand{\Supp}{\operatorname{Supp}}

\newcommand{\Pic}{\operatorname{Pic}\nolimits}
\newcommand{\pr}{\operatorname{pr}\nolimits}

\newcommand{\res}{\operatorname{res}}
\newcommand\Hilb{{\operatorname{Hilb}\nolimits}}

\newcommand{\CC}{{\mathbb C}}

\newcommand{\ZZ}{{\mathbb Z}}
\newcommand{\QQ}{{\mathbb Q}}
\newcommand{\PP}{{\mathbb P}}

\newcommand{\FFB}{{\mathbf F}}
\newcommand{\JJB}{{\mathbf J}}
\newcommand{\LLB}{{\mathbf L}}
\newcommand{\XXB}{{\mathbf X}}
\newcommand{\qB}{\boldsymbol{q}}
\newcommand{\pB}{\boldsymbol{p}}
\newcommand{\GB}{{\boldsymbol{\Gamma}}}

\newcommand{\OOO}{{\mathcal O}}

\newcommand{\III}{{\mathcal I}}
\newcommand{\GGG}{{\mathcal G}}
\newcommand{\EEE}{{\mathcal E}}
\newcommand{\HHH}{{\mathcal H}}
\newcommand{\LLL}{{\mathcal L}}
\newcommand{\FFF}{{\mathcal F}}

\newcommand{\NNN}{{\mathcal N}}
\newcommand{\MMM}{{\mathcal M}}
\newcommand{\RRR}{{\mathcal R}}
\newcommand{\SSS}{{\mathcal S}}
\newcommand{\TTT}{{\mathcal T}}

\newcommand{\UUU}{{\mathcal U}}

\newcommand{\tto}[1]{\xrightarrow{#1}}
\newcommand{\ffrom}[1]{\xleftarrow{#1}}

\newcommand{\isoto}{\stackrel{\textstyle\sim}{\lra}}

\newcommand{\cal}[1]{\mathcal{#1}}
\renewcommand{\bar}[1]{\overline{#1}}

\newcommand\alp{\alpha}

\newcommand\G{\Gamma}

\newcommand\si{\sigma}

\newcommand\lra{{\longrightarrow}}

\author{D. Markushevich}

\address{D. M.: Math\'ematiques - b\^{a}t. M2, Universit\'e Lille 1,
F-59655 Villeneuve d'Ascq Cedex, France}
\email{markushe@gat.univ-lille1.fr}

\author{A. S. Tikhomirov}

\thanks{Partially supported by the grant INTAS-OPEN-97-2072}

\address{A. T.: Department of Mathematics, State Pedagogical University,
Respublikanskaya 108, Yaroslavl' 150000, Russia}
\email{alexandr@tikho.yaroslavl.su}


\title{Symplectic structure on a moduli space of sheaves on
the cubic fourfold}

\begin{document}
\begin{abstract}
A 10-dimensional symplectic moduli space of torsion sheaves
on the cubic 4-fold is constructed. It parametrizes the stable
rank 2 vector bundles on the hypeplane sections of the cubic
4-fold which are obtained by Serre's construction from normal elliptic
quintics. The natural projection to the dual projective 5-space
parametrizing the hyperplane sections is a Lagrangian fibration.
The symplectic structure is closely related (and conjecturally,
is equal) to the quasi-symplectic one, induced by the
Yoneda pairing on the moduli space.
\end{abstract}
\maketitle

\section*{Introduction}

Mukai constructed in \cite{Muk} an example of a symplectic moduli space of
torsion sheaves on a K3 surface $S$, whose supports run
over a complete linear system on $S$. This moduli space is identic to
the relative compactified Jacobian of the linear
system, and the fibers, that is the Jacobians of the individual curves,
are Lagrangian subvarieties with respect to the symplectic structure.

We searched for an analogous situation with the relative intermediate
Jacobian of a family of higher dimensional varieties in place of the usual
Jacobian of a family of curves. Thus, the corresponding torsion sheaves
should be with supports of dimension at least 3, and their restrictions
to the supports should be of rank $>1$. We found a 
10-dimensional component of the moduli
space of stable sheaves whose supports are hyperplane sections
of a 4-dimensional cubic in $\PP^5$ with desired properties.
This is
a new example of a moduli space of sheaves possessing a symplectic
structure; the known ones parametrize sheaves on 
manifolds with holomorphic 2-forms (see \cite{Muk}
for K3 or abelien surfaces, \cite{K} for hyperkaehler varieties,
\cite{Tyu} for surfaces with $p_g>0$),
and the symplecitc structure comes essentially from the 2-form on the
base. In our example, the base
is very far from having holomorphic 2-forms: it is a Fano variety.

The fiber of our variety, parametrizing sheaves
with a fixed support (a cubic 3-fold), was constructed
in our previous paper \cite{Ma-Ti}. We studied there
the Abel--Jacobi map of the family of elliptic quintics lying
on a general cubic threefold $X$. We proved that
it factors through a 5-dimensional moduli 
component $M_X$ of stable rank 2 vector
bundles $\EEE$ on $X$ with Chern numbers $c_1=0,
c_2=2$, whose general point represents
a vector bundle obtained by Serre's construction from an elliptic
quintic $C\subset X$. The elliptic quintics mapped to a point of the moduli
space vary in a 5-dimensional projective space
inside the Hilbert scheme $\Hilb^{5n}_X$,
and the map $\Psi$ from the moduli space to the intermediate Jacobian $J(X)$
is \'etale on the open set representing (smooth) elliptic quintics which are
not contained in a hyperplane. This map is uniquely determined by the 
choice of a reference quintic in $X$. In fact, the elliptic quintics represent
$c_2(\EEE (1))$, so one can interprete $\Psi$ as $c_2$ with values
in algebraic cycles modulo rational equivalence, composed with the Abel--Jacobi
map on cycles.

In the present paper, we relativize this construction 
over the family of nonsingular
hyperplane sections of a cubic fourfold $Y$. This gives a 10-dimensional
moduli component $\frak M_Y$ of torsion sheaves on $Y$
with supports on the hypeplane sections $X=H\cap Y$, whose restrictions
to $X$ are rank two vector bundles from $M_X$. We show that
the Yoneda pairing induces a quasi-symplectic structure on it,
that is a holomorphic nondegenerate 2-form which is not
necessarily closed, and that
the natural projection $p:\frak M\lra\check{\PP}^5$
is a Lagrangian fibration, where
$\check{\PP}^5$ is the base of the family of hyperplanes in $\PP^5$.
This means that the fibers $M_X$, parametrizing sheaves
with fixed support $X$, are Lagrangian submanifolds of $\frak M$.

We study also the relation of the Yoneda quasi-symplectic structure
to the symplectic one 
constructed by Donagi--Markman \cite[Example 8.22]{D-M}
on $\JJB_Y$. It satisfies a similar property: the fibers 
$J(Y\cap H)$ are Lagrangian. In general,
there is no map between the two 10-dimensional varieties,
so we cannot hope that one quasi-symplectic structure is induced by the other.
Nevertheless, the map from $\frak M_Y$ to $\JJB_Y$ exists locally over
$\check{\PP}^5$
in the classical (or \'etale) topology, and we prove that these local
maps respect the symplectic structures up to the addition of a 2-form
which is a local section of $p^*\Omega^2_{\check{\PP}^5}$.
%This is the closest possible relation
%between two symplectic structures in the absence of a natural map from
%$\frak M_Y$ to $\JJB_Y$.

According to \cite{Il-Ma}, the degree of the above map $\Psi$ is
1, so $M_X$ can be identified with an open subset of $J^2(X)$,
and $\frak M_Y$ with that of some torsor under the relative
intermediate Jacobian $\JJB_Y$ of the
family of nonsingular hyperplane sections of $Y$.

Unfortunately, we do not see any reasonable approach to the proof of
the closedness of the Yoneda quasi-symplectic structure.
Remark, that Tyurin and Mukai did not prove the closedness of their
2-forms. Finally it turned out that their 2-forms are closed
\cite{O'G}, \cite{Bot}, \cite{Hu-Le}, but
the method for proving this in the cited papers
reduces  essentially the question to the closedness 
of the representatives for the Atiyah
class of a coherent sheaf. Our quasi-symplectic structure does not
involve components of the Atiyah class, so this method does not work.
Nevertheless, 
%after having
%established the relation to the symplectic structure of
%Donagi--Markman, 
we prove the existence of a symplectic
structure on $\frak M_Y$ which is in the
same relation to that of Donagi--Markman:
the above local maps between $\frak M_Y$ to $\JJB_Y$
respect the symplectic structures
up to the addition of a 2-form
which is a local section of $p^*\Omega^2_{\check{\PP}^5}$.
We define it in fact on the torsor, containing $\frak M_Y$.

Now, we will briefly describe the contents of the paper by sections.

In Section 1, we obtain the preparatory information on the cohomology
and the Ext groups of the sheaves from $\frak M_Y$. The calculations use
the local-to-global spectral sequence for Ext's and the projective
geometry of elliptic quintics.

In Section 2, we prove that the Yoneda pairing
$\Ext^1(\EEE ,\EEE )\times\Ext^1(\EEE ,\EEE )
\xrightarrow{\mbox{\scriptsize Yoneda}}
\Ext^2(\EEE ,\EEE )\simeq \CC$ for $\EEE\in{\frak M}_Y$ is skew-symmetric
and non-degenerate. We show that
the complex lines $\Ext^2(\EEE ,\EEE )$ for different $\EEE$'s fit together
into the trivial line bundle over ${\frak M}_Y$. Hence the Yoneda
coupling defines a genuine 2-form $\Lambda_{\frak M_Y}$ 
on ${\frak M}_Y$, which is the wanted
quasi-symplectic structure.

In Section 3, we describe the symplectic structure 
$\Lambda_{\JJB_Y}$ of Donagi--Markman
on $\JJB_Y$. It is determined by the one on the Fano 4-fold $Z$ of $Y$,
parametrizing lines $l\subset Y$. The existence and the uniqueness of the
latter was proved by
Beauville--Donagi \cite{B-D}, but they did not give any construction of it.
Section 4 provides an explicit construction
of the symplectic structure of Beauville--Donagi on $Z$. The idea
is to mimic the proof of the Tangent Bundle Theorem for the Fano
surface of a cubic 3-fold. Essentially, this proof is the calculation
of the differential of the Abel--Jacobi map of lines via the Clemens--
Griffiths--Welters diagram (see Sect. 2 in \cite{We}) for the triple
$l\subset X\subset\PP^4$. We adapt it to the triple $l\subset Y\subset\PP^5$,
and it gives (after an obvious relativization) a tensor field 
$\wedge_Z:\ \ \TTT Z\ \ {\tilde\to}\ \ \Omega_Z$.
The Bochner principle sais that the global holomorphic tensor fields
on $Z$ are generated by the Beauville--Donagi form, so $\wedge_Z$ has to
be proportional to it.

In Section 5, we compare the quasi-symplectic structures on $\frak M_Y$
and on $\JJB_Y$. These two manifolds are fibered over the same base $U
\subset\check{\PP}^5$ parametrizing the nonsingular hyperplane sections
of $Y$, have Lagrangian fibers,
and the fibers are related by the \'etale map
of \cite{Ma-Ti}.
So we can identify pointwise the tangent spaces of two respective fibers,
as well as their
normal spaces, which are identified with the tangent space to $U$. The first
thing to notice is that the two quasi-symplectic structures define the
same pairing between the tangent and the normal space of the fiber.
This gives a pointwise identification of two quasi-symplectic structures.
Next, we consider a local (in the classical topology)
fiberwise map $\Psi_\UUU : \frak M_Y
\to \JJB_Y$, which is the Abel--Jacobi map determined by a 
family of reference elliptic quintics in the hyperplane sections of $Y$ over
a small disk $\UUU\subset U$. We show that $\Lambda_{\frak M_Y}$,
$\Lambda_{\JJB_Y}$ differ by a section of $p^*\Omega^2_\UUU$.
In conclusion, we show that the Yoneda quasi-symplectic structure can be replaced by a symplectic structure with the same relation to $\Lambda_{\JJB_Y}$
and discuss the extensions to the bigger open set $\bar{U}$ parametrizing
hypeplane sections with at most one nondegenerate double point as
singularity.

\smallskip

{\em Acknowledgements}. The first author acknowledges with
pleasure the hospitality
of the ICTP at Trieste and of the MPIM at Bonn, where
he made a part of the work on the paper.

\section{Computing Ext's}

Let $Y\in\PP^5$ be a nonsingular cubic
fourfold, and $X=Y\cap\PP^4$ its general hyperplane section. Let $\EEE$
be the vector bundle obtained from a projectively normal elliptic quintic
curve $C\subset X\subset\PP^4$
by Serre's construction:
\begin{equation}\label{serre}
0\lra \OOO_X\lra \EEE (1) \lra \III_C(2) \lra 0\; ,
\end{equation}
where $\III_C=\III_{C,X}$ is the ideal sheaf of $C$ in $X$.
Since the class of $C$ modulo algebraic equivalence
is $5l$, where $l$ is the class of a line,
the sequence (\ref{serre}) implies that
$c_1(\EEE )=0, c_2(\EEE )=2l$. Moreover, $\det\EEE$ is trivial,
and hence  $\EEE$ is self-dual as soon as it is a vector
bundle (that is, $\EEE^\ast\simeq\EEE$).
See \cite[Sect. 2]{Ma-Ti} for further details on this construction.

Here, we will consider $\EEE$ as a torsion sheaf on $Y$, supported
on $X$.

\begin{lemma}\label{loc-ext}
 We have
$$
\EEE xt^q_{\OOO_Y}(\EEE ,\EEE )=\left\{
\begin{array}{l}
\EEE\otimes\EEE\ \mbox{\em if}\ q=0, \\
\EEE\otimes\EEE (1)\ \mbox{\em if}\ q=1, \\
0\ \ \mbox{\em otherwise}.
\end{array}\right.
$$
\end{lemma}

\begin{proof}
Choose a local trivialization of $\EEE$ :
$$
\EEE|_U\simeq (\OOO_X\oplus\OOO_X)|_U.
$$
This reduces the problem to the computaion of $\EEE xt^q_{\OOO_Y}
(\OOO_X ,\OOO_X)$. In using the resolution
$$
0\lra \OOO_Y(-1)\lra\OOO_Y\lra\OOO_X\lra 0,
$$
we arrive at the conclusion that $\EEE xt^q_{\OOO_Y}
(\OOO_X ,\OOO_X)=\OOO_X$ for $q=0$, $\NNN_{X/Y}\simeq\OOO_X(1)$ for $q=1$,
and $0$ for $q\geq 2$. Taking into accout the functorial
behavior of the last isomorphism with respect to the choice
of bases in $\OOO_X$, we obtain the wanted formulas
with $\EEE^*\otimes\EEE$ instead of $\EEE\otimes\EEE$, which
is the same by self-duality.
\end{proof}

\begin{lemma}\label{EE0}
The following assertions are true:

(i)
$
h^p(\EEE\otimes\EEE )=\dim\;\Ext^p_X(\EEE\otimes\EEE )=
\left\{
\begin{array}{l}
1\ \mbox{\em if}\ q=0, \\
5\ \mbox{\em if}\ q=1, \\
0\ \ \mbox{\em otherwise}.
\end{array}\right.
$

(ii) $h^0 (\EEE (1))=6$ and $h^i (\EEE (1))=0$ for $i>0$.

(iii) For $k\geq 1$,
$h^0(\III_C(k) )=
\binom{4+k}{4}-\binom{1+k}{4}
-5k$, and $h^i(\III_C(k) )=0$ if $i>0$.

(iv) $h^i(\EEE (j))=0$ for all $i\in\ZZ ,j=0,-1,-2$.

\end{lemma}

\begin{proof} (i),
(ii) follow from \cite[Lemma 2.1, c); Lemma 2.7, b)]{Ma-Ti}.
(iii) is verified for $k\geq 2$ in the proof of Lemma 2.1, loc. cit.
and is obviously extended to $k=1$ by considering the
exact triple
\begin{equation}\label{restriction}
0\lra\III_C(k)\lra\OOO_X(k)\lra\OOO_C(k)\lra 0
\end{equation}
with $k=1$. $C$ is not contained in a hyperplane by projective normality,
so $h^0 (\III_C(1))=0$. As $h^i (\OOO_X(1))=h^i(\OOO_C(1))$ equals
$5$ for $i=0$ and is zero for all $i>0$, we obtain the wanted assertion.

(iv) for $j=1$ is nothing but \cite[Lemma 2.7, a)]{Ma-Ti}.
By Serre duality (on $X$), it remains to verify (iv) only for $j=0$.
This follows from (iii) with $k=1$, the exact sequence (\ref{serre}),
twisted by $\OOO_X(-1)$, and from $H^i(\OOO_X(-1))=0$ for $i>0$.
\end{proof}

\begin{corollary}\label{E2pq}
The local-to-global spectral sequence
$$
E_2^{pq}=H^p(Y,\EEE xt^q_{\OOO_Y}(\EEE ,\EEE ))\Longrightarrow
\Ext^{p+q}_Y(\EEE ,\EEE )
$$
degenerates in its second term: $E_2^{pq}=E_{\infty}^{pq}$.
Hence $\Ext^0_Y(\EEE ,\EEE )=\CC$, that is, $\EEE$ is simple,
$\Ext^1_Y(\EEE ,\EEE )$ fits into the exact sequence
$$
0\lra H^1(\EEE\otimes\EEE )\lra \Ext^1_Y(\EEE ,\EEE )
\lra H^0(\EEE\otimes\EEE (1))\lra 0,
$$
and $\Ext^n_Y(\EEE ,\EEE )=H^{n-1}(\EEE\otimes\EEE (1))$ for
$n\geq 2$.
\end{corollary}

\begin{proof}
This follows immediately from Lemma \ref{loc-ext} and (i) of Lemma \ref{EE0}.
\end{proof}

\begin{lemma} \label{EE1}
The following assertions hold:

(i) $H^i(\EEE\otimes\EEE )=H^i(\EEE\otimes\III_C(1))$ and
$H^i(\EEE\otimes\EEE (1))=H^i(\EEE\otimes\III_C(2))$
for all $i\in\ZZ$.

(ii) $H^i(\EEE\otimes\EEE (1))=0$ for $i=2,3$, and $\chi (\EEE\otimes\EEE (1))=
h^0(\EEE\otimes\EEE (1))-h^1(\EEE\otimes\EEE (1))=4$.
%and $h^0(\EEE\otimes\EEE (1))\geq 4$.
\end{lemma}

\begin{proof} The first isomorphism is proved in \cite{Ma-Ti}.
For the second one, look at the exact sequences
\begin{equation}\label{EE}
0\lra \EEE \lra \EEE\otimes\EEE(1)\lra \III_C\otimes\EEE (2)\lra 0,
\end{equation}
\begin{equation}\label{restE}
0\lra \III_C\otimes\EEE (2)\lra \EEE (2)\lra \EEE (2)\otimes\OOO_C\lra 0,
\end{equation}
obtained by tensoring (\ref{serre}),(\ref{restriction}) with
$\EEE$.

The assertion (i) follows from (\ref{EE}) and Lemma \ref{EE0}, (iv).

%By Serre duality, $H^3(\EEE\otimes\EEE (1))=H^0(\EEE\otimes\EEE (-3))^*$.
%We have $H^0(\EEE\otimes\EEE )=
%\Ext^0_Y(\EEE ,\EEE )=\CC$, and hence $H^0(\EEE\otimes\EEE (-1))=0$;
%otherwise, $\dim H^0(\EEE\otimes\EEE )$ would be at least
%$h^0(\OOO_X (1))=5$. Hence $H^0(\EEE\otimes\EEE (-3))=0$.

The assertions (ii) and (iv) of Lemma \ref{EE0} give sufficiently many
values of $\chi (\EEE (k))$, which is a cubic polynomial
in $k$ by Riemann--Roch--Hirzebruch, that one can deduce that
$\chi (\EEE (k))=k(k+1)(k+2)$. Thus, $\chi (\EEE (2))=24$. We have also
$\deg\;\EEE (2)\otimes\OOO_C=20$, hence $\chi (\EEE (2)\otimes\OOO_C)=
20$, and this allows to determine the remaining Euler numbers from
(\ref{EE}),(\ref{restE}), thus proving $\chi (\EEE\otimes\EEE (1))=4$.

Now, use again (\ref{restE}) to show that $H^i(\EEE\otimes\III_C(2))=0$
for $i=2,3$. First, show that $H^1(\EEE (2)\otimes\OOO_C)=0$.
By (\ref{serre}), $\EEE (1)|_C\simeq\NNN_{C/X}$. Hence
$H^1(\EEE (2)\otimes\OOO_C)=H^1(\NNN_{C/X}(1))=H^0(\NNN_{C/X}^*(-1))^*=
H^0(\NNN_{C/X}(-3))^*$ (we used the fact that $\NNN_{C/X}(-1)$ is self-dual).
Consider the natural inclusion
$\NNN_{C/X}(-2)\subset\NNN_{C/\PP^4}(-2)$ and apply
\cite[Proposition V.2.1]{Hu}\label{V21}
: $H^1(\NNN^\ast_{C/\PP^4}(2)\otimes \MMM)=0$
for any invertible $\MMM$ of degree 0. By Serre duality,
$H^0(\NNN_{C/\PP^4}(-2))=0$.
Hence $H^0(\NNN_{C/X}(-2))=0$, and all the more $H^0(\NNN_{C/X}(-3))=0$.
Hence $H^1(\NNN_{C/X}(1))=0$.

Next, $H^i(\EEE (2))=0$ for $i>0$ by the same argument as the one
applied in the proof of (ii) of Lemma \ref{EE0}. Now (\ref{restE})
gives the wanted assertion.
\end{proof}

\begin{lemma}\label{IC2(3)}
We have, $h^0(\III_C^2(3))=0$, $h^1(\EEE (2)\otimes\III_C)=1$
and hence $h^0(\EEE (2)\otimes\III_C)=5$.
\end{lemma}

\begin{proof}
As in the proof of Lemma VI.3.1 in \cite{Hu}, we have the exact
sequence
\begin{equation}\label{hulek}
0\lra \FFF\lra H^0(\III_{C/\PP^4}(2))\otimes\OOO_C
\stackrel{\alpha}{\lra}\NNN_{C/\PP^4}^*(2)\lra 0,
\end{equation}
where $\FFF =\ker\:\alp$ is a rank two vector bundle
with $\det \FFF =\OOO_C(-1)$ and $h^0(\FFF )=0$. It is either the
direct sum of two line bundles $\FFF =\MMM_1\oplus\MMM_2$ with
$\deg M_i\geq -5$, or fits into the non-split exact sequence
$$
0\lra \OOO_C(-3o)\otimes\MMM\lra\FFF\lra\OOO_C(-2o)\otimes\MMM\lra 0,
$$
where $\MMM$ is a theta characteristic on $C$, and $o$ a point of
$C$ such that $\OOO_C(1)\simeq\OOO_C(5o)$.
This implies that $h^1(\FFF (1))=0$. By (\ref{hulek}), the map
$$
H^0(\III_{C/\PP^4}(2))\otimes H^0(\OOO_C (1))\lra H^0(\NNN_{C/\PP^4}^*(3))
$$
is surjective. This map factors obviously through
$H^0(\III_{C/\PP^4}(3))$, hence the map
$$
H^0(\III_{C/\PP^4}(3))
\lra H^0(\NNN_{C/\PP^4}^*(3))
$$
is also surjectve. By Riemann-Roch, $h^0(\FFF (1))=5$, and  by (\ref{hulek}),
$h^0(\NNN_{C/\PP^4}^*(3))=25-5=20$, and
$h^1(\NNN_{C/\PP^4}^*(3))=0$.
By \cite[Lemma 4.2]{Ma-Ti}, $h^0(\III_{C/\PP^4}(3))=20$,
hence, the exact sequence
$$
0\lra \III^2_{C/\PP^4}(3)\lra \III_{C/\PP^4}(3)\lra
\NNN_{C/\PP^4}^*(3)\lra 0
$$
yields $h^0(\III^2_{C/\PP^4}(3))=0$. By projective normality of $C$,
$h^1(\III_{C/\PP^4}(3))=0$, hence $h^1(\III^2_{C/\PP^4}(3))=0$.

Consider the following commutative diagram with exact rows and columns:
\begin{equation}\label{CD1}
\begin{CD}
 @. 0 @. 0 @. 0 @. \\
@. @VVV @VVV @VVV @. \\
0 @>>> \III_{C/\PP^4} @>>> \III^2_{C/\PP^4}(3) @>>>
\III^2_{C/X}(3) @>>> 0 \\
@. @VVV @VVV @VVV @. \\
0 @>>> \OOO_{\PP^4} @>>> \III_{C/\PP^4}(3) @>>> \III_{C/X}(3) @>>> 0 \\
@. @VVV @VVV @VVV @. \\
0 @>>> \OOO_C @>>> \NNN^*_{C/\PP^4}(3) @>>> \NNN^*_{C/X}(3) @>>> 0 \\
@. @VVV @VVV @VVV @. \\
 @. 0 @. 0 @. 0 @. \\
\end{CD}
\end{equation}

The left column gives $h^1 (\III_{C/\PP^4})=0$, and the first row
$h^0 (\III^2_{C/X}(3))= h^1 (\III_{C/\PP^4})=0$.
The lower row gives $h^0 (\NNN^*_{C/X}(3))=20$,
$h^1 (\NNN^*_{C/X}(3))=0$.
The middle row gives $h^0 (\III_{C/X}(3))=19$, $h^1 (\III_{C/X}(3))=0$,
which implies also $h^1 (\III^2_{C/X}(3))=1$.

Now, look at another commutative diagram in which all the triples
are exact:
\begin{equation}\label{CD2}
\begin{CD}
 @. @. 0 @. 0 @. \\
@. @. @VVV @VVV @. \\
0 @>>> \OOO_X(1) @>>> \EEE (2)\otimes\III_C @>>> \III_C^2(3) @>>> 0 \\
@. @| @VVV @VVV @. \\
0 @>>> \OOO_X(1) @>{\lambda}>> \EEE (2) @>{\mu}>> \III_C(3) @>>>  0 \\
@. @. @VVV @VVV @. \\
 @.  @. \NNN_{C/X}(1) @= \NNN^*_{C/X}(3) @.  \\
@. @. @VVV @VVV @. \\
 @. @. 0 @. 0 @. \\
\end{CD}
\end{equation}
From the middle row, we see that $\mu$ is
surjective on global sections.
We can think of elements of $H^0 (\EEE (2)\otimes\III_C)$
as of sections from $H^0 (\EEE (2))$ vanishing on $C$. Thus, a
section $s\in H^0 (\EEE (2))$ is in $H^0 (\EEE (2)\otimes\III_C)$
if and only if $\mu (s)\in H^0 (\III^2_{C/X}(3))$. Hence
$H^0 (\EEE (2)\otimes\III_C)=\mu^{-1}(H^0 (\III^2_{C/X}(3)))$
and
$$
 h^0 (\EEE (2)\otimes\III_C)=h^0 (\III^2_{C/X}(3)))+\dim\ker\mu =5,
$$
and
$h^1 (\EEE (2)\otimes\III_C)=1$ by Lemma \ref{EE1}, (ii).
\end{proof}

\begin{corollary}
$\dim \Ext^1_Y(\EEE ,\EEE )=10, \ \dim \Ext^2_Y(\EEE ,\EEE )=1$,
$\Ext^i_Y(\EEE ,\EEE )=0$ for $i>2$.
\end{corollary}

\begin{proof}
Corollary \ref{E2pq} and the last lemma imply the result.
\end{proof}

\section{Quasi-symplectic structure on $\frak M_Y$}

Let $Y$ be a nonsingular cubic fourfold. The notation $X$ will be used for
any of the nonsingular hyperplane sections
of $Y$. Denote by $\MMM$ the union of components of the moduli space of simple
sheaves on $Y$ containing the classes $[\EEE ]$ of the sheaves $\EEE$
defined by
(\ref{serre}). According to \cite{Ma-Ti}, the family
of elliptic normal quintics in a general cubic threefold
is irreducible, so $\MMM$ will be irreducible in the case if
$Y$ is general.
Let $\frak M\subset \MMM$ be the subset of classes of sheaves
$\EEE$ defined by (\ref{serre}) for all projectively normal elliptic
quintics $C$ lying in all nonsingular hyperplane sections of $Y$.

\begin{proposition}\label{moduli}
(i) $\frak M$ is an open subset of $\MMM$ contained in the smooth
locus of the moduli space of simple sheaves.

(ii) For all $i,p,q$, the groups $\Ext^i_Y(\EEE ,\EEE )$,
$H^p(Y,\EEE xt^q_{\OOO_Y}(\EEE ,\EEE ))$
for different $\EEE$ fit into
vector bundles $\TTT^i$, resp. $\TTT^{p,q}$ on $\frak M$.
That is, there exist such
vector bundles $\TTT^i$, $\TTT^{p,q}$ on $\frak M$ that for every $\EEE$
with $e=[\EEE ]\in\frak M$, there are canonical isomorphisms
$\TTT^i\otimes \CC (e)=\Ext^i_Y(\EEE ,\EEE )$,
$\TTT^{p,q}\otimes \CC (e)=H^p(Y,\EEE xt^q_{\OOO_Y}(\EEE ,\EEE ))$.
There is a natural including $\TTT^{i,0}\subset \TTT^i$.

$\TTT^1$ is canonically isomorphic to the tangent bundle of $\frak M$,
and $\rk\TTT^1=\dim\frak M = 10$, $\rk\TTT^2=1$.

(iii) The Yoneda coupling
$$
\Ext^i_Y(\EEE ,\EEE )\times\Ext^j_Y(\EEE ,\EEE )\lra
\Ext^{i+j}_Y(\EEE ,\EEE )
$$
defines a skew-symmetric bilinear form
\begin{equation}\label{skew}
\Lambda :\TTT^1\times\TTT^1\lra\TTT^2,
\end{equation}
and $\TTT^{1,0}$ is isotropic with respect to $\Lambda$.
\end{proposition}

\begin{proof}
(i) For the smoothness, remark that $\dim\frak M\leq 10=$ $
\dim \Ext^1_Y(\EEE ,\EEE )$, and there is an equality if and only
if $\MMM$ is smooth at $[\EEE ]$. But it is indeed an equality,
because $\frak M$ contains explicitly a 10-dimensional family of
pairwise non-isomorphic sheaves: we have 5 parameters specifying the
support $X$ of the sheaf, and the moduli
space $M_X$ with fixed support $X$ is smooth of dimension 5 by
\cite{Ma-Ti}. The openness of $\frak M$ in the moduli space of
sheaves on $Y$ follows exactly as in loc. cit., Corollary 5.5
in using an appropriate relativization of Serre's construction.

The existence of the sheaves $\TTT^i$ in (ii) and
of the bilinear form in (iii) follow from
Proposition 2.2 of \cite{Muk} and from Theorem 1.21, (4) of \cite{Sim}.
The existence of the $\TTT^{p,q}$ is obtained by an obvious modification
of the arguments of \cite{Muk}.
The skew symmetry follows from the smoothness:
the Yoneda square $z\circ z$ of an element $z\in\Ext^1_Y(\EEE ,\EEE )$
is exactly the first obstruction map $\Ext^1_Y(\EEE ,\EEE )\lra$  $
\Ext^2_Y(\EEE ,\EEE )$, and it is zero by smoothness of $\MMM$ at $\EEE$.

By general properties of spectral sequences of composite
functors with a ring structure, the Yoneda coupling
induces on the bi-graded object $E_2^{pq}$ a natural ring structure
$$
H^p(Y,\EEE xt^q_{\OOO_Y}(\EEE ,\EEE ))\times
H^{p'}(Y,\EEE xt^{q'}_{\OOO_Y}(\EEE ,\EEE ))\lra
H^{p+p'}(Y,\EEE xt^{q+q'}_{\OOO_Y}(\EEE ,\EEE ))
$$
The isotropy of $H^1(\EEE\otimes\EEE )$ follows from the fact that
$H^2(\EEE\otimes\EEE )=0$ (Lemma \ref{EE0}, (i)).
\end{proof}

In order to get a quasi-symplectic
structure, we have to verify
that $\Lambda$ is  non-degenerate, and that $\TTT^2$ is a trivial
line bundle.

%\subsection{Non-degeneracy of the coupling}

By the isotropy of $H^1(\EEE\otimes\EEE )$, the non-degeneracy
of the coupling $\Ext^1_Y(\EEE ,\EEE )\times\Ext^1_Y(\EEE ,\EEE )\lra
\Ext^2_Y(\EEE ,\EEE )$ is equivalent to that of
\begin{equation}\label{coupling0}
H^1(\EEE\otimes\EEE )\times H^0(\EEE\otimes\EEE (1))
\lra H^1(\EEE\otimes\EEE (1))
\end{equation}

Now, in using the isomorphism $\wedge^2\EEE\simeq\OOO_X$,
we can identify $H^0(\EEE\otimes\EEE (1))$ with
$H^0(\wedge^2\EEE (1))=H^0(\OOO_X(1))$ via the natural map
$\EEE\otimes\EEE \lra \wedge^2\EEE$. The fact that the resulting
map is an isomorphism follows from the factorisation
$$
H^0(\EEE\otimes\EEE (1))=H^0(\III_C\otimes\EEE (2))=H^0(\OOO_X(1)),
$$
coming from (\ref{EE}) and from the upper row
of (\ref{CD2}).

\begin{lemma}\label{plain-mult}
The coupling (\ref{coupling0}) is identified via the above
isomorphism with the coupling of plain multiplication of elements
of $H^1(\EEE\otimes\EEE )$ by linear forms:
\begin{equation}\label{coupling}
H^1(\EEE\otimes\EEE )\times H^0(\OOO_X(1))\lra H^1(\EEE\otimes\EEE (1)).
\end{equation}
\end{lemma}

\begin{proof}
To see that these are indeed the same, we will express both
of them in a local basis of $\EEE$. The latter is just the multiplication
of coefficients:
$$
(\sum a_{ij} e_i\otimes e_j\ ,\ l
\ )\ \mapsto\ \ \sum la_{ij} e_i\otimes e_j\ .
$$
For the former, notice that $\EEE\otimes\EEE$ contains the trivial
subbundle $\wedge^2\EEE$ as a direct summand. By the above
considerations, $H^0(\EEE\otimes\EEE (1))=H^0(\wedge^2\EEE (1))$. Hence
every section from $H^0(\EEE\otimes\EEE (1))$ can be written locally
as $l(e_2\otimes e_1-e_1\otimes e_2)$. The self-duality of $\EEE$ being given
by a trivialization of $\wedge^2\EEE$, we can assume that $e_2\wedge e_1
=e_2\otimes e_1-e_1\otimes e_2$
is the trivializing section, and then the dual basis of
$(e_j)=(e_1,e_2)$ is $((-1)^{\bar{j}}e_{\bar{j}})=(e_2,-e_1)$,
where $\bar{j}=3-j$. Hence the contraction of middle terms
$$
(a_1\otimes a_2, b_1\otimes b_2)\mapsto \frac{a_2\wedge b_1}{e_2\wedge e_1}
a_1\otimes b_2
$$
giving the Yoneda coupling is the identity $e_i\otimes e_j
\mapsto e_i\otimes e_j$ on the basic sections.
This implies that the Yoneda coupling is given by multiplication
of coefficients, $(a_{ij}, l)\mapsto la_{ij}$.
\end{proof}

\begin{proposition}
The coupling (\ref{coupling}) is non-degenerate.
\end{proposition}

\begin{proof}
We are to verify, that for any linear form $l\in H^0(\OOO_X(1))$,
the map
$$
\cdot l :H^1(\EEE\otimes\EEE )\lra H^1(\EEE\otimes\EEE (1))(\simeq \CC )
$$
is surjective.
%Consider the exact triple
%\begin{equation}\label{cdot-l}
%0\lra \III_C\otimes\EEE (1)\stackrel{\cdot l}{\lra} \III_C\otimes\EEE (2)
%\lra\III_C\otimes\EEE (2)|_H\lra 0,
%\end{equation}
%where the first map is the multiplication by a linear form
%$l$, and $H$ is the hyperplane section defined by $l$.
%Assume that
%$l$ is sufficiently general, so that $H$ is a smooth cubic
%surface and $S:=H\cap C$ is a set of five distinct points.
%Write $\EEE_H=\EEE |_H$ and $\III_S=\III_{S/H}=\III_C|_H$.
By Lemma 1.4, the last map is identified with
$$
\cdot l :H^1(\III_C\otimes\EEE (1))\lra H^1(\III_C\otimes\EEE (2)) .
$$

Using (\ref{CD1}), as it is and twisted by $\OOO (-1)$, and
Lemma~\ref{EE1}, (i),
we can continue the sequence of identifications:
\begin{equation}\label{CD5terms}
\begin{CD}
H^1(\III_C\otimes\EEE (1)) @>{\sim}>>
H^1 (\III^2_{C/X}(2)) @>{\sim}>>
H^2 (\III_{C/\PP^4}(-1)) @<{\sim}<<
H^1(\OOO_C(-1)) \\
@V{\cdot l}VV @V{\cdot l}VV @V{\cdot l}VV @V{\cdot l}VV      \\
H^1(\III_C\otimes\EEE (2)) @>{\sim}>>
H^1 (\III^2_{C/X}(3)) @>{\sim}>>
H^2 (\III_{C/\PP^4})
@<{\sim}<<
H^1(\OOO_C ) \\
\end{CD}
\end{equation}
Only the isomorphisms involving $\III^2_{C/X}$ are not obvious.
Those involving $H^1 (\III^2_{C/X}(3))$ are established in the
proof of Lemma 1.5. We will verify now the isomorphisms involving
$H^1 (\III^2_{C/X}(2))$.
To this end, write down the diagram (\ref{CD1}) twisted by $\OOO (-1)$:

\begin{equation}\label{CD1'}
\begin{CD}
 @. 0 @. 0 @. 0 @. \\
@. @VVV @VVV @VVV @. \\
0 @>>> \III_{C/\PP^4}(-1) @>>> \III^2_{C/\PP^4}(2) @>>>
\III^2_{C/X}(2) @>>> 0 \\
@. @VVV @VVV @VVV @. \\
0 @>>> \OOO_{\PP^4}(-1) @>>> \III_{C/\PP^4}(2) @>>> \III_{C/X}(2) @>>> 0 \\
@. @VVV @VVV @VVV @. \\
0 @>>> \OOO_C(-1) @>>> \NNN^*_{C/\PP^4}(2) @>>> \NNN^*_{C/X}(2) @>>> 0 \\
@. @VVV @VVV @VVV @. \\
 @. 0 @. 0 @. 0 @. \\
\end{CD}
\end{equation}

The left column yields for the vector $(h^i(\III_{C/\PP^4}(-1)))$
the values $(0,0,5,0,\ldots )$. Similarly the left column of (\ref{CD1})
gives $(h^i(\III_{C/\PP^4}(2)))=(5,0,0,\ldots )$. Now, the middle
row of (\ref{CD1'}) gives $(h^i(\III_{C/X}(2)))=(5,0,$\linebreak $0,\ldots )$.
Further, $\NNN^*_{C/X}(2)\simeq \NNN_{C/X}$, 
and by \cite[Lemma 5.1, e)]{Ma-Ti},
$(h^i(\NNN^*_{C/X}(2)))=(10,0,0,\ldots )$. We proved in the previous
section that $h^0(\III^2_{C/X}(3))=0$, and so much the more
$h^0(\III^2_{C/X}(2))=0$. Hence, from the last column of (\ref{CD1'})
we can deduce that $(h^i(\III^2_{C/X}(2)))=(0,5,0,\ldots )$.

We see already that $h^2(\III_{C/\PP^4}(-1))=h^1(\III^2_{C/X}(2))=5$.
Now, let us verify that the connecting homomorphism $H^1\lra H^2$ of the first
row is indeed an isomorphism. Look at the second column of (\ref{CD1'}).
By \cite[Proposition V.2.1]{Hu},
$(h^i(\NNN^*_{C/\PP^4}(2)))=(5,0,0,\ldots )$. We saw in the previous
section that $h^0(\III^2_{C/\PP^4}(3))=0$, and so much the more
$h^0(\III^2_{C/\PP^4}(2))=0$. Hence
$(h^i(\III^2_{C/\PP^4}(2)))=(0,0,0,\ldots )$. Hence all the
connecting homomorphisms of the first row are isomorphisms.

Now, coming back to (\ref{CD5terms}), we see that the vertical arrow
on the right hand side is obviously surjective, and we are done.
\end{proof}

\begin{theorem}
The line bundle $\TTT^2$ on $\frak M$ is trivial, hence the
skew-symmetric pairing $\Lambda$ defined in (\ref{skew}) is a
well defined $2$-form on $\frak M$:
\begin{equation}
\Lambda:\TTT{\frak M}\overset{\sim}{\to}\Omega^1_{\frak M},\ \ \
\Lambda\in H^0(\Omega^2_{\frak M}).
\end{equation}
The fibers of the natural map $p:\frak M\to\check\PP^5$, 
$\EEE\to <\Supp (\EEE)>$, are Lagrangian subvarieties with respect to $\Lambda$.

\end{theorem}

\begin{proof}
By Proposition \ref{moduli}, for any point
$e=[\EEE]\in\frak M$,
one has a canonical isomorphism
$\TTT^2\otimes{\bf C}(e)\simeq \Ext^2_Y(\EEE,\EEE)$.
It is enough to show that there is a canonical isomorphism
\begin{equation}\label{canon-iso}
\alpha:\ \Ext^2_Y(\EEE,\EEE)\overset{can}{\simeq}{\bf C}.
\end{equation}
To this end, fix an equation
$f_Y=0$
of the cubic $Y$ in $\PP^5$
and take any section
$h\in H^0(\OOO_Y(1))$
vanishing on
$X=\Supp(\EEE)$;
this section $h$ defines an isomorphism
$\phi_h:\OOO_X(1)\overset{\sim}{\to}\NNN_{X/Y}$
such that the substitution
$(h\mapsto\lambda h)$
induces
$\phi_h\mapsto\lambda^{-1}\phi_h$.

Next, by Lemma \ref{loc-ext} there is a canonical isomorphism
\begin{equation}\label{canon1}
\EEE xt^1_{\OOO_Y}(\EEE,\EEE)\overset{can}{\xrightarrow{\sim}}
\EEE^\vee\otimes\EEE\otimes\NNN_{X/Y}.
\end{equation}

Now take a general section
$0\ne s\in H^0(\EEE(1))$
such that
$C=(s)_0$
is a smooth projectively normal elliptic quintic.
By computations of Section 1,
\begin{equation}\label{vanish}
h^i(\det\EEE\otimes\III^2_{C,H}(3))=h^i(\III^2_{C,H}(3))=0,\ \ i=1,2,
\end{equation}
where
$H=<C>=<X>$ is the span of $C$, or of $X$, in $\PP^5$.
The section $s$ defines an exact sequence
\begin{equation}\label{seq1}
0\to\OOO_X\overset{s}{\to}\EEE(1)\overset{\cdot\wedge s}{\to}
\III_{C,X}(2)\otimes\det\EEE\to0.
\end{equation}
Twisting it by
$\EEE^\vee$
and using the canonical isomorphism
$\EEE^\vee\otimes\det\EEE\overset{can}{\xrightarrow{\sim}}\EEE$,
we obtain the exact triple
$0\to\EEE^\vee\to\EEE^\vee\otimes\EEE(1)
\overset{\varepsilon_s}{\to}\EEE\otimes\III_{C/X}(2)\to0$.
Thus we  have a composed isomorphism

\begin{multline}\label{compos1}
\gamma_1: \Ext^2_Y(\EEE,\EEE)\overset{can}{\underset{Cor. \ref{E2pq}}{\simeq}}
H^1(\EEE xt^1_{\OOO_Y}(\EEE,\EEE))\overset{can}{\simeq}
H^1(\EEE^\vee\otimes\EEE\otimes\NNN_{X/Y})
\underset{\sim}{\xrightarrow{\phi_h^{-1}}} \\
H^1(\EEE^\vee\otimes\EEE(1))
\underset{\sim}{\xrightarrow{H^1(\varepsilon_s)}}
H^1(\EEE\otimes\III_{C,X}(2))
\end{multline}

which is proportional, by the above, to
$h$ and to $s$
\footnote{We denote this as follows: $\gamma_1\sim(h,s)$.}.
Next, similar to (\ref{seq1}) we have a morphism
$\EEE\otimes\III_{C,X}(2)\overset{\cdot\wedge s}{\to}
\III^2_{C,X}(3)\otimes\det\EEE$
inducing a map
$H^1(\EEE\otimes\III_{C,X}(2))\to H^1(\III^2_{C,X}(3)\otimes\det\EEE)$
which gives, in composition with $\gamma_1$, the isomorphism

\begin{equation}\label{compos2}
\gamma_2: \Ext^2_Y(\EEE,\EEE)\overset{\sim}{\to}
H^1(\III^2_{C,X}(3)\otimes\det\EEE),\ \ \ \gamma_2\sim(h,s^2).
\end{equation}

The exact triples
$0\to\III_{C,H}\overset{\cdot f_Y}{\to}\III^2_{C,H}(3)\to\III^2_{C,X}(3)\to0$
and
$0\to\III_{C,H}\to\OOO_H\to\OOO_C\to0$
give together with (\ref{vanish}) and Serre duality on
$C$
a chain of canonical isomorphisms:

\begin{equation}\label{compos3}
H^1(\III^2_{C,X}(3)\otimes\det\EEE)\simeq
H^2(\III_{C,H}\otimes\det\EEE)\simeq
H^1(\det\EEE\otimes\OOO_C)\simeq
H^0((\det\EEE)^\vee\otimes\omega_C)^\vee.
\end{equation}

Next, 
$\phi_h:\OOO_X(1)\overset{\sim}{\to}\NNN_{X/Y},$
$\phi_h\sim h$,
composed with 
$\psi:\NNN_{C/X}\overset{\sim}{\to}\EEE(1)|_C,\ \ \psi\sim s,$
and with the canonical isomorphisms
$\omega_C\simeq\omega_Y\otimes\NNN_{X/Y}\otimes\det\NNN_{C/X},\ \ \
\omega_Y\simeq\OOO_Y(-3)$, gives the following isomorphism:
$\chi:\omega_C\overset{\sim}{\to}\det\EEE|_C,\ \ \ \ \chi\sim(h,s^2).$
Composing it with (\ref{compos3}), we obtain the isomorphism

\begin{multline}\label{compos4}
\theta:\ H^1(\III^2_{C,X}(3)\otimes\det\EEE)\overset{can}{\simeq}
H^0((\det\EEE)^\vee\otimes\omega_C)^\vee\overset{\sim}{\to}
H^0((\det\EEE|_C)^\vee\otimes\det\EEE|_C)\overset{can}{\simeq}{\bf C},\\
\ \ \ \ \theta\sim(h^{-1},s^{-2}).
\end{multline}

Now, composing (\ref{compos2}) and (\ref{compos4}), we get the desired
isomorphism $\alpha$ in (\ref{canon-iso}) with
$\alpha\sim(h^0,s^0),$
that is
$\alpha$
does not depend on the choice of
scalar multiples of $h$ and $s$. Thus, to show that
$\alpha$
is a canonical isomorphism, it is enough to check that
$\alpha$
does not depend on the choice of the point
$z={\bf C}s\in \PP (H^0(\EEE(1)).$
The last assertion is clear, since, by construction,
$\alpha$ is defined for {\em any} point
$z\in \PP (H^0(\EEE(1))$
\footnote
{In fact, for any
$z={\bf C}s\in \PP (H^0(\EEE(1))$, it follows
from the stability and from the local freeness of $\EEE$, that the scheme
$C=(s)_0$ of zeros of $s$
is a locally complete intersection of pure codimension 2 in $X$,
i.e. a curve, so that the map
$\alpha:\ \Ext^2_Y(\EEE,\EEE)\to{\bf C}$
is well defined.},
thus giving a morphism of trivial bundles
$\alpha:\ \Ext^2_Y(\EEE,\EEE)\otimes\OOO_{\PP (H^0(\EEE(1))}{\to}
{\bf C}\otimes\OOO_{\PP (H^0(\EEE(1))}$.
We have
$\alpha\in H^0(\OOO_{\PP (H^0(\EEE(1))})$,
hence $\alpha$ is constant as a function on the projective 5-space
$\PP (H^0(\EEE(1))$.

The tangent distribution of the fibers of $p$ is the vector
bundle $\TTT^{1,0}$, and its isotropy was proved in Proposition \ref{moduli}.
This implies that the fibers of $p$ are Lagrangian.

\end{proof}

\section{Symplectic structure on $\JJB_Y$}\label{DMMT}

In this section, $Y$ will denote a nonsingular cubic fourfold
in $\PP^5$,
$h:\JJB =\JJB_Y\lra U$ the relative intermediate Jacobian of the family
$\XXB\lra U\subset \check{\PP}^5$ of nonsingular hyperplane sections
over an open subset $U$ in the dual projective space $\check{\PP}^5$.
$X$ will denote any of the smooth hyperplane sections $H\cap Y$ for
$H\in U$.
$\XXB$ is the scheme
of zeros of a section of $\OOO_\Pi (3,0)$ in
the incidence divisor $\Pi =\{(x,H)\in\PP^5_U\mid x\in H\}$, where
$\PP^5_U=\PP^5\times U\subset\PP^5\times\check{\PP}^5$. We can fix this section,
as well as the cubic form defining $Y$ in $\PP^5$,
and consider all the isomorphisms depending only on these data
as canonical ones.

According to Donagi--Markman \cite[8.22]{D-M}, $\JJB$ has a natural
(up to a constant factor) symplectic
structure $\alp_\JJB\in H^0(\JJB ,\Omega^2_\JJB )$ for which $h$ is
a Lagrangian fibration, that is the fibers $J_X$ of $h$,
the intermediate Jacobians of the hyperplane sections $X=H\cap Y$,
are Lagrangian submanifolds of $\JJB$. The idea is to recover the
symplectic structure from the isomorphisms $\NNN_{L/\JJB}\simeq
\Omega_L^1$ that it defines for Lagrangian submanifolds $L$. For any symplectic
structure $\alp$ on $\JJB$ for which $h$ is a Lagrangian fibration,
the above isomorphisms for Lagrangian submanifolds $L=J_X$ glue
together to give an isomorphism of vector bundles
$\sigma (\alp ):\TTT_U\lra h_*\Omega^1_{\JJB /U}$.

Conversely, let $\sigma\in \Homg (\TTT_U,h_*\Omega^1_{\JJB /U})=
H^0(U,\Omega^1_U\otimes h_*\Omega^1_{\JJB /U})$ be any isomorphism.
The natural exact triple
\begin{equation}
0\lra h^*\Omega^1_U\lra \Omega^1_\JJB\lra\Omega^1_{\JJB /U}\lra 0
\end{equation}
implies another one, obtained by taking $\wedge^2$:
\begin{equation}\label{omega2U}
0\lra h^*\Omega^2_U\lra \FFF\lra h^*\Omega^1_U\otimes \Omega^1_{\JJB /U}\lra 0,
\end{equation}
where $\FFF \subset\Omega^2_\JJB$ is the vector subbundle of sections
for which $\TTT_{\JJB /U}$ is an isotropic distribution.
Let $\delta : H^0( \JJB , h^*\Omega^1_U\otimes \Omega^1_{\JJB /U})
\lra H^1(\JJB , h^*\Omega^2_U)$ be the connecting homomorphism defined
by (\ref{omega2U}). The following assertion is obvious:

\begin{lemma}\label{lem-obs}
$h^*\sigma$ lifts to a section $\alp\in H^0(\JJB ,\FFF )$ if and only
if $\delta (h^*\sigma )=0$. $\alp$ is determined up to the
addition of a 2-form $h^*\beta$, where $\beta\in H^0(U ,\Omega^2_U)$.
The lift becomes unique, if one imposes the condition that the
zero section $O\subset\/B$ be a Lagrangian subvariety (that is,
$\alp |_O\equiv 0$).
\end{lemma}

We will formulate the result of Donagi--Markman in a form, convenient
to our purposes. Let $Z$ be the ``Fano scheme" of $Y$, that is the
Hilbert scheme parametrizing all the lines in $Y$. It is a smooth projective
irreducible 4-fold. Beauville--Donagi \cite{B-D}
proved that $Z$ is an irreducible
symplectic variety, deformation equivalent to the Fujiki--Beauville
4-fold $S^{[2]}$ for a K3 surface $S$; let $\alp_Z$ be its symplectic
structure, unique up to proportionality. Let $\pi :\FFB\lra U$
be the relative Fano surface of $\XXB$ over $U$, that is the
family of surfaces $F(X)$ parametrizing the lines in the hyperplane sections
$X=H\cap Y$ with $H\in U$. By Clemens--Griffiths \cite{CG},
$\JJB /U$ is canonically isomorphic to the relative Albanese variety
$\Alb (\FFB /U)$ and to $\Pic^0(\FFB /U)$. 
By Voisin \cite{V}, the Fano surfaces $F(X)$
are Lagrangian in $Z$ and form an open set in the Hilbert scheme
of $Z$. By Theorem 8.1 of \cite{D-M}, the relative Picard variety
$P=\Pic (\LLB /B)$ of the universal family of Lagrangian subvarieties
of a smooth projective symplectic variety $V$ over an open set
$B\subset\Hilb (V)$ has also a symplectic structure for which
$P/B$ is a Lagrangian fibration. The application of this result to
the family $\FFB /U$ yields a symplectic 10-fold $\Pic^0(\FFB /U)=
\Alb (\FFB /U)=\JJB$.
\smallskip

{\bf Theorem of Donagi--Markman}. For any $H\in U$, let $X=Y\cap H$,
$F=F(X)$, $\sigma_H :\NNN_{F/Z}\lra\Omega^1_F$ the isomorphism,
induced by $\alp_Z$, and $H^0(\sigma_H ):H^0(\NNN_{F/Z})\lra H^0(\Omega^1_F)$
the induced map on global sections. With natural identificatons
$H^0(F,\NNN_{F/Z})=T_{[F]}\Hilb (Z)=T_{[F]}U$ and
$H^0(F,\Omega^1_F)=\Omega^1_{\Alb (F)}|_0=\Omega^1_{J(X)}|_0$, let $\sigma :\TTT_U\lra
h_*\Omega^1_{\JJB /U}$ be the isomorphism obtained by the relativization
of the construction of $H^0(\sigma_H )$ over $U$. Then $\sigma$ satisfies
the hypotheses of Lemma \ref{lem-obs}, and hence lifts to
a unique symplectic structure $\alp_\JJB$ on $\JJB$ for which the
fibers and the zero section of $h$ are Lagrangian submanifolds.\smallskip

Our goal now is to identify $\sigma$ as the canonical isomorphism,
constructed in Lemma \ref{canis} below.

Let $\pr_1,\pr_2$ be the natural projections onto the two factors of
$\PP^5\times U$, and $\OOO (i,j)=\pr_1^*\OOO_{\PP^5}(i)
\otimes\pr_2^*\OOO_{\check{\PP}^5}(j)$.

\begin{lemma}\label{canis}
There exist canonical isomorphisms
$$
h_*\Omega^1_{\JJB /U}=\pi_*\Omega^1_{\FFB /U}=
R^1\pr_{2*}\Omega^2_{\XXB /U}=\pr_{2*} \OOO_{\XXB}(1,1)=\TTT_U .
$$
\end{lemma}

\begin{proof} The first two isomorphisms follow from the definitions
of the intermediate Jacobian and of the Albanese variety, as well
as the identification $\JJB =\Alb (\FFB /U)$ mentioned above.
The third one is a  relative version of the tangent bundle formula
for the intermediate Jacobian of a cubic 3-fold. To prove it, write down
the relative conormal bundle sequence of $\XXB\subset\Pi$ over
$U$, taken to the third exteriour power and twisted by $\NNN_{\XXB /\Pi}$:
$$
0\lra \Omega^2_{\XXB /U}\lra\Omega^3_{\Pi /U}\otimes\NNN_{\XXB /\Pi}\lra
\Omega^3_{\XXB /U}\otimes\NNN_{\XXB /\Pi}\lra 0 .
$$
Applying $R^i\pr_{2*}$ and using the Bott formulae, we obtain
the canonical isomorphism
\begin{equation}\label{R1}
\pr_{2*}(\Omega^3_{\XXB /U}\otimes\NNN_{\XXB /\Pi})\isoto
R^1\pr_{2*}\Omega^2_{\XXB /U} .
\end{equation}
We have $\NNN_{\XXB /\Pi}=\OOO_\XXB (3,0)=\OOO_\Pi (3,0)|_{\XXB}$,
 $\Omega^3_{\XXB /U}=\omega_{\Pi /U}\otimes\NNN_{\XXB /\Pi}$
and $\omega_{\Pi /U}=\OOO_\Pi (-5,1)$, so both sheaves in (\ref{R1})
are canonically isomorphic to $\pr_{2*} \OOO_{\XXB}(1,1)$.

It is also well known that the incidence divisor $\Pi$ is identified
with the projectivized cotangent bundle of $U$:
$\Pi =\PP (\TTT_U^\dual)$, hence there are canonical isomorphisms
$\TTT^\dual_U=\pr_{2*}\OOO_{\Pi /U}(1)=\pr_{2*} \OOO_{\Pi}(1,1)
=\pr_{2*} \OOO_{\XXB}(1,1)$.

\end{proof}

The Clemens--Griffiths Tangent Bundle Theorem for the Fano surface
$F(X)$ of a cubic 3-fold $X$ states, that there is an isomorphism
$\TTT_{F(X)}\simeq\SSS_{G(1,\PP^4)}|_{F(X)}$, where $G(1,\PP^n)$
denotes the Grassmannian of lines in $\PP^n$ and $\SSS_{G(1,\PP^n)}$
the universal rank 2 vector bundle on it.
We are going to relativize this result, in order to compute the
sheaf $\Omega^1_{\FFB /U}$.

Let $G=G(1,\PP^5)$, $\SSS_G$ the universal vector bundle,
$\PP^5\ffrom{q_G}\G_G\tto{p_G}G$ the universal family of lines
together with its natural projections,
$Y\ffrom{q}\G\tto{p}Z$ its restriction to the lines in $Y$,
and $\XXB\ffrom{\qB}\GB\tto{\pB}\FFB$ the relativization over $U$.
Let $p_i$ be the projection to the three factors of
$\PP^5\times \check{\PP}^5\times G$, and $\OOO(i,j,k)=
p_1^*\OOO (i)\otimes p_2^*\OOO (j)\otimes p_3^*\OOO (k)$.
Denote by $\G_u$, $F_u$, $X_u$, $H_u=\PP^4_u$ the
fibers over $u\in U$ of $\GB$, $\FFB$, $\XXB$, resp. $\Pi$.

\begin{lemma}\label{FFB1}
There is a canonical isomorphism
$$
\TTT_{\FFB /U}=(\OOO_{\check{\PP}^5}(-1)\boxtimes\SSS_G)|_\FFB ,
$$
where the operation $\boxtimes$ yields a sheaf on the variety $\check{\PP}^5
\times G$, containing $\FFB$ in a natural way.
\end{lemma}

\begin{proof}
The natural inclusions
$\GB\subset \XXB\times_U\FFB\subset \Pi\times_U \FFB$
give the following exact sequences:
\begin{equation}\label{chert}
0\lra \NNN_{\GB /\XXB\times_U\FFB}\lra
\NNN_{\GB /\Pi\times_U\FFB}\lra
\NNN_{\XXB\times_U\FFB /\Pi\times_U\FFB}|_{\GB}\lra 0 ,
\end{equation}
$$
0\lra \TTT_{\GB /\FFB}\lra
(p_1\times p_2)^*\TTT_{\Pi /U}|_{\GB}\lra \NNN_{\GB /\Pi\times_U\FFB}\lra 0 .
$$
We have $\TTT_{\GB /\FFB}=(p_1\times p_3)^*\TTT_{\G_G/G}$. The Euler
sequence
$$
0\lra \OOO\lra q_G^*\OOO_{\PP^5}(1)\otimes p_G^*\SSS_G\lra
\TTT_{\G_G/G}\lra 0
$$
implies $\det \TTT_{\G_G/G}=q_G^*\OOO_{\PP^5}(1)\otimes p_G^*\OOO_G(-1)$.
We find successively the determinants:
$\det \TTT_{\GB /\FFB}=\OOO_{\GB} (2,0,-1)$,
$\det (p_1\times p_2)^*\TTT_{\Pi /U}|_{\GB}=\OOO_{\GB} (5,-1,0)$,
$\det \NNN_{\GB /\Pi\times_U\FFB}=\OOO_{\GB} (3,-1,1)$,
$\NNN_{\XXB\times_U\FFB /\Pi\times_U\FFB}|_{\GB}=
p_1^*\NNN_{Y/\PP^5}|_{\GB}=\OOO_{\GB} (3,0,0)$,
$\det \NNN_{\GB /\XXB\times_U\FFB}=\OOO_{\GB} (0,-1,1)$.
Dualizing (\ref{chert}) and tensoring by $\det \NNN_{\GB /\XXB\times_U\FFB}$,
we obtain the exact triple
$$
0\lra \OOO_{\GB} (-3,-1,1)\lra \NNN^{\dual}_{\GB /\Pi\times_U\FFB}(0,-1,-1)\lra
\NNN_{\GB /\XXB\times_U\FFB}\lra 0,
$$
which gives the canonical isomorphism
$$
\pB_*\NNN_{\GB /\XXB\times_U\FFB}\isoto R^1\pB_*\OOO_{\GB} (-3,-1,1) .
$$
As $\pB_*\NNN_{\GB /\XXB\times_U\FFB}=\TTT_{\FFB /U}$ and by
relaive Serre duality, we obtain the canonical isomorphism
$$
\TTT_{\FFB /U}=(\pB_*\OOO_{\GB} (1,1,0))^{\dual}
. $$
The natural identification $\SSS_G=(p_{G*}q_G^*\OOO_{\PP^5}(1))^\dual$
induces the canonical isomorphism
$(\pB_*\OOO_{\GB} (1,1,0))^*=(\OOO_{\check{\PP}^5}(-1)\boxtimes\SSS_G)|_\FFB$,
which implies the assertion of the lemma.

\end{proof}

The next lemma relativizes the isomorphism
$\NNN_{F_u/Z}\simeq\SSS_G^\dual|_{F_u}$
coming from $\NNN_{G(1,\PP^4_u)/G(1,\PP^5)}=\SSS_G^\dual|_{G(1,\PP^4_u)}$.

\begin{lemma}\label{canis2}
There exist canonical isomorphisms
$$
\NNN_{\FFB /\check{\PP}^5\times Z}=\OOO_{\check{\PP}^5}(1)
\boxtimes\SSS_G^\dual|_\FFB =\Omega^1_{\FFB /U} .
$$
Moreover, applying $\pi_*$, one obtains the isomorphism
$$
\pi_*\NNN_{\FFB /\check{\PP}^5\times Z}=
\pi_*\Omega^1_{\FFB /U} ,
$$
which coincides with that of Lemma \ref{canis} if one
uses the natural identification
$\pi_*\NNN_{\FFB /\check{\PP}^5\times Z}=\TTT_U$,
following from the interpretation of $\FFB /U$ as the
universal family over the open set $U\subset\Hilb (Z)$.
\end{lemma}

\begin{proof}
One can represent $\FFB$ as the complete intersection
$\boldsymbol{\sigma_{11}}\cap (U\times Z)$, where
$$
\boldsymbol{\sigma_{11}}=\{ (u,l)\in \check{\PP}^5\times G
\mid l\subset H_u\}
$$
is the universal Schubert variety $\si_{11}$ over $\check{\PP}^5$.
It is the scheme of zeros of $(q_G\times\id )_*(p_G\times\id )^* \tau_0$,
where $\tau_0\in H^0(\PP^5\times\check{\PP}^5, \OOO_{\PP^5\times\check{\PP}^5}
(1,1))$ is the equation of $\Pi$. As
$(q_G\times\id )_*(p_G\times\id )^*\OOO_{\PP^5\times\check{\PP}^5}
(1,1)=\OOO_{\check{\PP}^5}(1)\boxtimes\SSS_G^\dual$,
we have proved the first statement.

The second one follows immediately from the Clemens--Welters
computation of the differential of the Abel--Jacobi map (see
Sect. 2 in \cite{We}, or diagram (\ref{CDA8}) in Sect. \ref{section5} below).
\end{proof}

%Using the identifications of Lemma \ref{canis2},
%Write down the cotangent bundle sequence for the embedding $\FFB\subset
%U\times Z$ over $U$:
%\begin{equation} \label{CBS2}
%0\lra \NNN^\dual_{\FFB /U\times Z}
%\lra \ka^*\Om^1_Z\lra\Om^1_{\FFB /U}\lra 0,
%\end{equation}
%where $\ka :\FFB\lra Z$ is the natural projection.
%Applying to (\ref{CBS2}) $\wedge^2$, one gets the following exact
%triple, similar to (\ref{omega2U}):
%\begin{equation}\label{lambda2T}
%0\lra\wedge^2\NNN^\dual_{\FFB /U\times Z}\lra \GGG \lra
%\NNN^\dual_{\FFB /U\times Z}\otimes\Om^1_{\FFB /U}\lra 0 ,
%\end{equation}
%where $\GGG $ is a vector subbundle of $\ka^*\Om^2_Z$.

\section{Symplectic structure on the Fano variety $Z$ of
the cubic 4-fold $Y$}

\subsection{Construction of the symplectic structure}

Let $\overline{\mathbf F}$, resp. $\overline{\boldsymbol{\Gamma}}$
be the closure of
$\FFB$
in
$\check\PP^5\times G$,
resp.
of
$\GB$
in
$\PP^5\times\check\PP^5\times G$.
\renewcommand{\FFB}{\overline{\mathbf F}}
\renewcommand{\GB}{\overline{\boldsymbol{\Gamma}}}

We have the diagram of natural embeddings

\begin{equation}\label{CDA1}
\begin{CD}
@. Y\times\FFB @.\ \ \subset\ \  @. \PP^5\times\FFB @. \\
@. \cup @. @. \cup @. \\
\GB\ \ \subset\ \  @. \XXB\times_{\check\PP^5}\FFB @.
\ \ \subset\ \  @. \Pi\times_{\check\PP^5}\FFB @. \\
\end{CD}
\end{equation}

Remark, that $\FFB$ and $\GB$ are nonsingular of dimensions
7 and 8 respectively. This follows from the fact that $\FFB$
admits a natural smooth morphism onto $Z$ with fibers $\PP^3$,
and $\GB /\FFB$ is smooth with fibers $\PP^1$. Write down the
exact sequence of normal bundles:

\begin{equation}\label{CDA2}
\begin{CD}
0 @>>> \NNN_{\GB /Y\times\FFB} @>>>
\NNN_{\GB /\PP^5\times\FFB} @>>>
\NNN_{Y\times\FFB/\PP^5\times\FFB}|_{\GB} @>>> 0 
\end{CD}
\end{equation}

Applying $\wedge^2$ and using the isomorphism
$
\NNN_{Y\times\FFB/\PP^5\times\FFB}|_{\GB}
\simeq\OOO_{\GB}(3,0,0),
$
we obtain the exact triple

\begin{equation}\label{CDA3}
\begin{CD}
0 @>>> \wedge^2\NNN_{\GB /Y\times\FFB}
@>>> \wedge^2\NNN_{\GB /\PP^5\times\FFB}
@>>> \NNN_{\GB /Y\times\FFB}(3,0,0)
@>>> 0
\end{CD}
\end{equation}

Now, from the proof of Lemma \ref{FFB1}, it follows
$
\det\NNN_{\GB /{\check\PP^5}\times\FFB}=
\det \NNN_{\GB /\XXB\times_{\check\PP^5}\FFB}\otimes
\det\NNN_{Y\times\FFB/\PP^5\times\FFB}|_{\GB}=\OOO_{\GB} (1,0,1),
\omega_{\GB /\FFB}=\OOO_{\GB} (-2,0,1).
$
Using these formulas and twisting the last exact triple by
$\OOO_{\GB}(-3,0,0)$, we obtain:

\begin{equation}\label{CDA4}
\begin{CD}
0 @>>> \NNN^\vee_{\GB /Y\times\FFB}\otimes\omega_{\GB /\FFB}
@>>> \wedge^2\NNN_{\GB /\PP^5\times\FFB}(-3,0,0)
@>>> \NNN_{\GB /Y\times\FFB}
@>>> 0
\end{CD}
\end{equation}

\begin{remark}\label{remA}
Consider the commutative square

\begin{equation}\label{CDA5}
\begin{CD}
\GB @.\ \ \overset{\tilde\kappa}{\lra}\ \  @. \Gamma \\
\ \ \downarrow\pB @. @. \ \ \downarrow p \\
\FFB @. \ \ \overset{\kappa}{\lra}\ \  @. Z, \\
\end{CD}
\end{equation}

where $\tilde\kappa$ is the restriction to $\GB$ of the natural projection
$\PP^5\times\check\PP^5\times G\to\PP^5\times G$.
Then (\ref{CDA4})
%\begin{equation}\label{middle}
%0\to\NNN^\vee_{\GB /Y\times\FFB}\otimes\omega_{\GB /\FFB}
%\to\wedge^2\NNN_{\GB /\PP^5\times\FFB}(-3,0,0)
%\to\NNN_{\GB /Y\times\FFB}
%\to0
%\end{equation}
is obtained by applying
$\tilde\kappa^*$
to the triple

\begin{equation}\label{mid}
0\to\NNN^\vee_{\Gamma/Y\times Z}\otimes\omega_{\Gamma/Z}
\to\wedge^2\NNN_{\Gamma/\PP^5\times Z}(-3,0)
\to\NNN_{\Gamma/Y\times Z}
\to0.
\end{equation}

\end{remark}

For any point
$z\in\FFB$,
respectively,
$l\in Z$,
we have
$
\wedge^2\NNN_{\GB /\PP^5\times\FFB}(-3,0,0)|\pB^{-1}(z)
\simeq\OOO_{\PP^1}(-1)^{\oplus6},\ \
$,
respectively,
$
\wedge^2\NNN_{\Gamma/\PP^5\times Z}(-3,0)|p^{-1}(l)
\simeq\OOO_{\PP^1}(-1)^{\oplus6},
$
so that the base change implies
\begin{equation}\label{R^i}
R^i{\pB}_*(\wedge^2\NNN_{\GB /\PP^5\times\FFB}(-3,0,0))=0,\ \
R^ip_*(\wedge^2\NNN_{\Gamma/\PP^5\times Z}(-3,0))=0,\ \ i\ge0,
\end{equation}

Applying
$R^ip_*$ to (\ref{mid}) and,
respectively,
$R^i{\pB}_*$
to (\ref{CDA4}), and using
(\ref{R^i}), relative Serre duality, Remark \ref{remA}
and the isomorphism
$p_*\NNN_{\Gamma/\PP^5\times Z}\simeq\TTT Z$,
we obtain the isomorphisms

\begin{equation}\label{sympl}
\wedge_Z:\ \ \TTT Z\ \ {\tilde\to}\ \ \Omega_Z,
\end{equation}

and, respectively,
\begin{equation}\label{symplect}
\kappa^*\wedge_Z:\ \ \kappa^*\TTT Z\ \ {\tilde\to}\ \ \kappa^*\Omega_Z.
\end{equation}

Since, by \cite{B-D} $Z$ is an irreducible symplectic variety,
the Bochner principle for irreducible symplectic manifolds
\cite{B} implies that (\ref{sympl})
is a symplectic form on $Z$. Thus we
obtain an explicit construction of the symplectic structure on $Z$,
which is stated in the
following theorem:

\begin{theorem}\label{symplform}
The symplectic structure on $Z$ is the connecting homomorphism (\ref{sympl})
in the long exact sequence of the functors $R^ip_*$, associated
to the exact triple (\ref{mid}).
\end{theorem}

\renewcommand{\FFB}{{\mathbf F}}
\renewcommand{\GB}{{\boldsymbol{\Gamma}}}

Now recall that, by the results of the previous section,
\begin{equation}\label{pB_*}
\pB_*\NNN_{\GB /\XXB\times_{ U}\FFB}=\TTT\FFB/U,
\ \ \pB_*\OOO_\GB(1,1,0)=\NNN_{\FFB/U\times Z},
\end{equation}
By relative Serre duality,
\begin{equation}\label{R^1pB_*}
R^1\pB_*(\NNN^\vee_{\GB /\XXB\times_{U}\FFB}
\otimes\omega_{\GB /\FFB})=\Omega_{\FFB/U},
\ \ R^1\pB_*(\OOO_\GB(1,1,0)^\vee\otimes\omega_{\GB /\FFB})=
\NNN^\vee_{\FFB/U\times Z}.
\end{equation}
Thus, using (\ref{pB_*}), (\ref{R^1pB_*}), (\ref{symplect})
and applying $R^i\pB_*$ to (\ref{CDA4}),
we obtain the commutative diagram:

\begin{equation}\label{CDA6}
\begin{CD}
0 @>>> \TTT\FFB/U
@>>> \kappa^*\TTT Z
@>>> \NNN_{\FFB/U\times Z} @>>> 0 \\
@. @VVV \kappa^*\wedge_Z@VVV \sigma_0@VVV @. \\
0 @>>> \NNN^\vee_{\FFB/U\times Z}
@>>> \kappa^*\Omega_Z
@>>> \Omega_{\FFB/U}
@>>> 0 \\
\end{CD}
\end{equation}

By Theorem \ref{symplform}, $\wedge_Z$ is a symplectic form;
besides, by construction, $\sigma_0$ is the isomorphism of Lemma
\ref{canis2}. Thus we obtain

\begin{theorem}\label{ASD}
The commutative diagram (\ref{CDA6}) is antiselfdual and the isomorphism
$\sigma_0$ in this diagram coincides with the canonical isomorphism of
Lemma \ref{canis2}.
\end{theorem}

In particular, the composition
$
\TTT\FFB/U \to \kappa^*\TTT Z
\ \overset{\kappa^*\wedge_Z}{\lra}\
\kappa^*\Omega_Z \to \Omega_{\FFB/U}
$
in this diagram is the zero map. This gives another proof of the result of
Voisin \cite{V} that Fano surfaces of hyperplane sections of $Y$
are the Lagrangian subvarieties of the symplectic structure on $Z$.

\section{Relation between symplectic structures on $\frak M_Y$ and $\JJB_Y$}
\label{section5}
\subsection{Codifferential of the Abel-Jacobi map.} We recall here the
following Clemens-Welters interpretation of the codifferential
$d\Phi^\vee$
of the Abel-Jacobi map
$\Phi: B\ \to J(X)$,
where $B$ is the base of a certain family of curves on a smooth hyperplane
section $X=Y\cap\PP^4$ of the cubic $Y$ (see \cite[Sect. 2]{We}).
Let $[C]$ be a point in
$B$
corresponding to a locally complete intersection (l.c.i.) curve
$C\subset X$
and let $y=\Phi(C)$. Applying $\wedge^3$ to the exact triple
$
0\to\NNN^\vee_{X/\PP^4}\to\Omega^1_{\PP^4}|X\to\Omega^1_X\to0
$
and twisting by
$\NNN_{X/\PP^4}=\OOO_X(3)$,
we obtain the exact triple
\begin{equation}\label{wedge3X}
0\to\Omega^2_X\to\Omega^3_{\PP^4}\otimes\NNN_{X/\PP^4}\to\OOO_X(1)\to0.
\end{equation}
Passing to the cohomology of this sequence, we obtain the coboundary isomorphism
\begin{equation}\label{iso-R}
R:\ \ H^0(\OOO_X(1))\ \ \overset{\sim}{\to}\ \ H^1(\Omega^2_X).
\end{equation}
Remark that, by construction, $R$ is the restriction to a fiber
${\PP^4}$
of the isomorphism
$\pr_{2*} \OOO_{\XXB}(1,1)=R^1\pr_{2*}\Omega^2_{\XXB /U}$
from Lemma \ref{canis}.

Next, the triple (\ref{wedge3X}) restricted to $C$ clearly
fits into the commutative diagram

\begin{equation}\label{CDA7}
\begin{CD}
0 @>>> \Omega^2_X|C
@>>> \Omega^3_{\PP^4}\otimes\NNN_{X/\PP^4}|C
@>>> \OOO_C(1)
@>>> 0 \\
@. @VVV @VVV \parallel @. @. \\
0 @>>> \Omega^3_X\otimes N_{C/X}
@>>> \Omega^3_{\PP^4}\otimes N_{C/\PP^4}
@>>> \Omega^3_X\otimes N_{X/\PP^4}|C
@>>> 0 \\
\end{CD}
\end{equation}

Remark that, for a l.c.i. curve $C\subset X$,
$\Omega^3_X\otimes N_{C/X}=\omega_X\otimes N_{C/X}=
\omega_X\otimes \det N_{C/X}\otimes N^\vee_{C/X}=
\omega_C\otimes N^\vee_{C/X}$, so that by Serre duality
$H^1(\Omega^3_X\otimes N_{C/X})=H^0(N_{C/X})^\vee$.

Hence, passing to the cohomology in the bottom triple,
we obtain the map

\begin{equation}\label{beta_C}
\beta_C:\ \ H^0(\OOO_C(1))\ \ \to\ \ H^0(N_{C/X})^\vee.
\end{equation}

According to loc. cit., the codifferential
$d\Phi^\vee|_{[C]}$
at the point $[C]$ of the Abel-Jacobi map
$\Phi: B\ \to\ J(X)$
makes the following diagram commutative:

\begin{equation}\label{CDA8}
\xymatrix{
H^0(\OOO_X(1)) \ar[d]^{\otimes\OOO_C} \ar[rr]^R_\sim & &
H^1(\Omega^2_X) \ar[d]^\wr \\
H^0(\OOO_C(1)) \ar[r]^{\beta_C} & H^0(\NNN_{C/X})^\vee\ \  &
\Omega^1_{J(X)}|y \ar[l]_{d\Phi^{\vee}|\{ C\}}.
}
\end{equation}

Now, as in Section \ref{DMMT}, let
$U=\{\PP^4\in\check\PP^5|X=Y\cap\PP^4$ is smooth$\}$,
$h:\JJB\to U$
and
$\pi:\FFB\to U$
the projections and
$\sigma=\pi_*\sigma_0:\pi_*\NNN_{\FFB/\check\PP^5\times Z}|U
\overset{\sim}{\to}\pi_*\Omega^1_{\FFB/\check\PP^5}|U
\simeq h_*\Omega^1_{\JJB/U}$
the canonical isomorphism (see Lemmas \ref{canis} and \ref{canis2}).
Take a point
$y\in h^{-1}(U)$ and let $X=\PP^4\cap Y$,
where
$\PP^4=h(y).$
Remark that, by Lemma \ref{canis2}, there is a natural identification
$H^0(\OOO_X(1))\simeq h^*\pi_*\NNN_{\FFB/\check\PP^5\times Z}|y$.
These isomorphisms together with (\ref{CDA8}) and the corresponding base
change isomorphisms give the commutative diagram
\begin{equation}\label{CDA9}
\xymatrix{
H^0(\OOO_X(1)) \ar[d]^\wr_R \ar[rrrr]^\simeq & & & &
h^*\pi_*\NNN_{\FFB/\check\PP^5\times Z}|y \ar[d]^{h^*\sigma}_\wr \\
H^1(\Omega^2_X) & \simeq & \Omega^1_{\JJB/U}|y &
\simeq & h^*h_*\Omega^1_{\JJB/U}|y \\
}
\end{equation}
Here $H^0(\OOO_X(1))=h^*\TTT U|y$ and the isomorphism $R$ globalizes to
\begin{equation}\label{RRR}
\RRR:\ \ h^*\TTT U\ \ \overset{\sim}{\to}\ \ \Omega^1_{\JJB/U}.
\end{equation}

Now consider the symplectic structure of Donagi-Markman
$\wedge_\JJB:\TTT\JJB\overset{\sim}{\to}\Omega^1_{\JJB}$
induced on $\JJB$ via the symplectic structure $\wedge_Z$ on $Z$.
Then putting together (\ref{CDA9}) and the right hand commutative square of 
(\ref{CDA6}), we obtain

\begin{theorem}\label{sympl-R}
There is an antiselfdual commutative diagram
\begin{equation}\label{CDA10}
\begin{CD}
0 @>>> \TTT\JJB/U
@>>> \TTT\JJB
@>>> h^*\TTT U @>>> 0 \\
@. |\wr@. \wedge_\JJB\ |\wr@.  \RRR\ |\wr@. @. \\
0 @>>> h^*\Omega^1_{U}
@>>> \Omega^1_{\JJB}
@>>> \Omega^1_{\JJB/U}
@>>> 0 \\
\end{CD}
\end{equation}
\end{theorem}

\subsection{Local lifts of $\wedge_\JJB$ to $\frak M$}

Consider  the open subset
$
\HHH'=\{c\in \Hilb^{5n}_Y|\ \ C$
is a smooth projectively normal quintic
$\}$
of $\Hilb^{5n}_Y
$,
with its projection
$\rho':\HHH'\to\check\PP^5:C\mapsto <C>$.
Let $Y$ be general enough, so that, by \cite{Ma-Ti},
$\HHH'$ is an irreducible smooth open subset of
$\Hilb^{5n}_Y, \ \ \dim\HHH'=15,\ \ \rho(\HHH')$
is dense open in $\check\PP^5$ and there exists
a dense open subset
$U\subset\check\PP^5$,
respectively, a dense open subset
$\HHH$ of $\HHH'$, $\HHH\subset{\rho'}^{-1}(U)$,
such that any $X\in U$ is a smooth hyperplane section of $Y$
and
$\rho=\rho'|\HHH:\HHH\to U$
is a smooth morphism of relative dimension 10.
For any point
$C\in\HHH$,
there exists a local analytic (or local in the \'etale topology) section of the projection
$\rho$ passing through $C$, i.e. an analytic  (resp. \'etale) neighbourhood
$\UUU\ni\rho(C)$
and a holomorphic (resp. regular) map
$s=s_\UUU :\UUU\to\HHH$
such that
$s(\UUU)\ni C$.
Let
$
\HHH_{\UUU}=\rho^{-1}(\UUU),\ \
\JJB_{\UUU}=h^{-1}(\UUU),\ \ \frak M_{\UUU}=
p^{-1}(\UUU),
$
where
$p:\frak M\to\check\PP^5:\EEE\to <\Supp (\EEE)>$
is the projection, and denote
$p_{\UUU}=p|_{\frak M_{\UUU}}$.
The section $s$ defines a  relative Abel-Jacobi map
$\Phi_s:\HHH_{\UUU}\to\JJB_{\UUU}:C\mapsto\Phi_{X,s}(C)$,
where
$\Phi_{X,s}:\rho^{-1}(X)\to J_X:C\mapsto\Phi(C-s(X))$
is the usual Abel-Jacobi map,
$J_X=h^{-1}(X)$.
By \cite[Theorem 5.6]{Ma-Ti} the map
$\Phi_\UUU$ has the Stein factorization
$\Phi_\UUU=\Psi_\UUU\circ\phi_\UUU$,
where
$\phi_\UUU:\HHH_\UUU\to\frak M_\UUU$
is a smooth holomorphic map of relative dimension 5 and
$\Psi_\UUU:\frak M_\UUU\to\JJB_\UUU$
is an \'etale holomorphic map.

Next, passing to $\wedge^2$ in the lower triple of (\ref{CDA10}),
we obtain the diagram

\begin{equation}\label{CDA11}
\begin{CD}
@. 0 @. @. @. \\
@. @VVV @. @. \\
@. h^*\Omega^2_\UUU @. @. @. \\
@. @VVV @. @. \\
0 @>>> \FFF
@>>> \Omega^2_{\JJB}
@>{e_\JJB}>> \Omega^2_{\JJB/\UUU}
@>>> 0 \\
@. @V{\varepsilon_\JJB}VV @. @. \\
@. \HHH om(h^*\TTT\UUU,\Omega^1_{\JJB_\UUU/\UUU})@. @. @. \\
@. @VVV @. @. \\
@. 0 @. @. @. \\
\end{CD}
\end{equation}

Similarly we have on ${\frak M}_\UUU$

\begin{equation}\label{CDA12}
\begin{CD}
@. 0 @. @. @. \\
@. @VVV @. @. \\
@. p_\UUU^*\Omega^2_\UUU @. @. @. \\
@. @VVV @. @. \\
0 @>>> \GGG
@>>> \Omega^2_{{\frak M}_\UUU}
@>{e_{\frak M}}>> \Omega^2_{{\frak M}_\UUU/\UUU}
@>>> 0, \\
@. @V{\varepsilon}VV @. @. \\
@. \HHH om(p_\UUU^*\TTT\UUU,\Omega^1_{{\frak M}_\UUU})@. @. @. \\
@. @VVV @. @. \\
@. 0 @. @. @. \\
\end{CD}
\end{equation}
and since
$\Psi_\UUU:\frak M_\UUU\to\JJB_\UUU$
is \'etale, the diagram (\ref{CDA12}) can be obtained from (\ref{CDA11})
by applying
$\Psi_\UUU^*$:
\begin{equation}\label{relation}
(\ref{CDA12})=\Psi_\UUU^*(\ref{CDA11}).
\end{equation}

Now let
$\Lambda_\JJB\in H^0(\Omega^2_\JJB)$
be the Donagi-Markman symplectic structure on $\JJB$
and
$\Lambda_{\JJB_\UUU}=\Lambda_\JJB|\frak M_{\UUU}$
its restriction to
$\JJB_{\UUU}$,
and let $\Lambda_{\frak M_{\UUU}}=\Lambda|\frak M_{\UUU}$,
where
$\Lambda$ is the Yoneda quasi-symplectic form introduced
in Proposition \ref{moduli} above).
By Proposition \ref{moduli}, for any
$[\EEE]\in{\frak M}_\UUU$
the vector space
$T_{[\EEE]}{\frak M}_\UUU/\UUU\subset T_{[\EEE]}{\frak M}_\UUU$
is isotropic w.r.t.
$\Lambda_{{\frak M}_\UUU}$,
hence

\begin{equation}\label{isotr1}
e_{\frak M}(\Lambda_{{\frak M}_\UUU})=0.
\end{equation}

By the same reason
\begin{equation}\label{isotr2}
e_\JJB(\Lambda_{\JJB_\UUU})=0.
\end{equation}
The conditions (\ref{isotr1}) and (\ref{isotr2}) imply that
$\Lambda_{{\frak M}_\UUU}\in H^0(\GGG )$,
respectively,
$\Lambda_{\JJB_\UUU}\in H^0(\FFF )$.
Now prove that
\begin{equation}\label{equality1}
\varepsilon(\Lambda_{{\frak M}_\UUU})=
\Psi_\UUU^*\varepsilon_\JJB(\Lambda_{\JJB_\UUU}),
\end{equation}
or equivalently, according to (\ref{relation}), that
\begin{equation}\label{equality2}
\varepsilon(\Lambda_{{\frak M}_\UUU}-\Psi_\UUU^*\Lambda_{\JJB_\UUU})=0.
\end{equation}

To prove (\ref{equality1}), take any
$\EEE\in{\frak M}_\UUU$
and let
$y=\Psi_\UUU([\EEE]),\ \ X=h(y).$
By virtue of the identifications
$T_{[\EEE]}{\frak M}_\UUU/\UUU=H^1(\EEE\otimes\EEE),\ \
p_\UUU^*\TTT\UUU|_{[\EEE]}=H^0(\OOO_X(1))$
we have an isomorphism
${}^\sharp\Lambda:=\varepsilon(\Lambda_{\frak M_\UUU})|_{[\EEE]}:
H^0(\OOO_X(1))\overset{\sim}{\to}H^1(\EEE\otimes\EEE)^\vee$.
By Theorem \ref{sympl-R} together with the identifications,
$T_y\JJB_\UUU/\UUU=H^1(\Omega^2_X)^\vee,\ \
h^*\TTT\UUU|_y=H^0(\OOO_X(1))$
we have an isomorphism
$R=\varepsilon_\JJB(\Lambda_{\JJB_\UUU})|_y:
H^0(\OOO_X(1))\overset{\sim}{\to}H^1(\Omega^2_X)$.
Thus, by (\ref{CDA8}), the equality (\ref{equality1}) is reduced to the
commutativity of the diagram

\begin{equation}
\xymatrix{
H^1(\Omega^2_X) \ar[r]^{d\Psi^\dual_\UUU} & H^1(\EEE\otimes\EEE)^\vee \\
H^0(\OOO_X(1))  \ar[u]_\wr^R \ar[ur]_{{}^\sharp\Lambda}&\\
} \ \ \ ,
\end{equation}

i.e., after dualizing and using(\ref{CDA8}),
to that of

\begin{equation}
\xymatrix{
H^0(\NNN_{C/X}) \ar[d]_{\beta_C^\vee} \ar[r]^{d\Phi_{X,s}|_{[C]}} &
H^1(\EEE\otimes\EEE) \ar[d]^{({}^\sharp\Lambda)^\vee} \\
H^0(\OOO_C(1))^\vee \ar[r]^\sim &
H^0(\OOO_X(1))^\vee\\
} \ \ \ ,
\end{equation}
where
$C$ is any curve in
$\phi^{-1}_\UUU([\EEE])$.

By Lemma \ref{plain-mult} and (\ref{CD5terms}),
the last diagram is equivalent to

\begin{equation}\label{CDA13}
\xymatrix{
H^0(\NNN_{C/X})\otimes H^0(\OOO_C(1))
\ar[d]_{d\Phi_{X,s}|_{[C]}\otimes(\res^{-1})} \ar[r]^{\ \ \ \ \ m} &
H^0(\NNN_{C/X}(1)) \ar[d]^{\tilde\delta} \\
H^1(\EEE\otimes\EEE)\otimes H^0(\OOO_X1)) \ar[r]^{\ \ \ \ \ m} &
H^1(\EEE\otimes\EEE(1)),\\
}
\end{equation}
where $m$ is the ordinary multiplication,
$\tilde\delta$
the composition
$
\tilde\delta:\ \ H^0(\NNN_{C/X}(1))\overset{\delta}{\twoheadrightarrow}
H^1(\III_{C,X}^2(3))\overset{\gamma_1}{\underset{\sim}{\longleftarrow}}
H^1(\III_{C,X}(2)\otimes\EEE)\underset{\sim}{\overset{\gamma_2}{\longleftarrow}}
H^1(\EEE\otimes\EEE(1)),
$\ \ \ \
$\delta$ the connecting map in the cohomology of the exact triple
$0\to\III_{C,X}^2(3))\to\III_{C,X}(3))\to\NNN_{C/X}(1)\to0$,
and the isomorphisms
$\gamma_1,\ \gamma_2$
are obtained from the triples
$0\to\OOO_X(1)\to\III_{C,X}(2)\otimes\EEE\to\III_{C,X}^2(3)\to0,\ \
0\to\EEE\to\EEE\otimes\EEE(1)\to\III_{C,X}(2)\otimes\EEE\to0$
in using the equalities
$h^i(\OOO_X(1))=0,\ \ h^i(\EEE)=0,\ i=1,2$. 
The commutativity of (\ref{CDA13}) is obvious.
Hence, (\ref{equality1}) follows.

Denote
$\Lambda_{\UUU,s}=\Lambda_{{\frak M}_\UUU}-\Psi_\UUU^*\Lambda_{\JJB_\UUU}.$
From (\ref{equality2}) and (\ref{CDA11}), it follows that
$\Lambda_{\UUU,s}\in H^0(p_\UUU^*\Omega^2_\UUU).$
\footnote{Clearly, the form
$\Lambda_{\UUU,s}$
depends on the choice of the section $s:\UUU\to\HHH.$}
Thus, we obtain the following theorem relating
the quasi-symplectic structures
$\Lambda_{{\frak M}_\UUU}$ and $\Lambda_{\JJB_\UUU}$.

\begin{theorem}\label{relat}
For the data $(\UUU,s)$ specified above, there exists a 2-form
$\Lambda_{\UUU,s}\in H^0(p_\UUU^*\Omega^2_\UUU)$
such that
\begin{equation}
\Lambda_{{\frak M}_\UUU}=\Psi_\UUU^*\Lambda_{\JJB_\UUU}+\Lambda_{\UUU,s}
\end{equation}
\end{theorem}

\subsection{Symplectic structure on ${\frak M}$}

Only at this point we will use the result of \cite{Il-Ma}:

\begin{theorem}
For any nonsingular cubic threefold $X$, the map $\Psi :M_X\lra J^2(X)$ is
an open embedding.
\end{theorem}

\begin{proof}
See Theorem 3.2 and Corollaries 3.3, 5.1 in \cite{Il-Ma}.
\end{proof}

This theorem implies that the above maps $\Psi_\UUU :{\mathfrak M}_\UUU\lra\JJB_\UUU$
defined by local sections $s_\UUU :\UUU\lra \HHH$ of $\rho :\HHH\lra U$
are also open embeddings. Let $\{\UUU_\alp\}_{\alp\in I}$ be
an open covering of $U$
in the classiacal or \'etale topology such that there exist local sections $s_\alp :\UUU_\alp\lra\HHH$
of $\rho$, and let $\si_\alp :\UUU_\alp\lra
{\mathfrak M}_\UUU$ be the corresponding sections of $\Psi$. Then $c=(\si_{\alp\beta})_{\alp ,\beta\in I}$,
where $\si_{\alp\beta}=\si_{\beta}-\si_{\alp}$, is a 1-cocycle of the covering $\{\UUU_\alp\}$ with values in
the sheaf $\underline{\JJB}$ of holomorphic (or regular) sections of $\JJB$, and
${\mathfrak M}$ is naturally identified with an open subset in the torsor
$\JJB^c$ constructed from $c$. The following lemma 
summarizes some easy facts about the relations between 
structures of Lagrangian fibrations
on an abelian scheme and on its torsors.

\begin{lemma}\label{torsors}
Let $B$ be a nonsingular variety of dimension $n$ and $f:A\lra B$ a smooth
family of connected commutative algebraic (or complex Lie) groups such that
the generic fiber $A_b$ is an abelian variety (or a compact complex torus) of dimension $n$. Let $c$ be any cohomology class from
$H^1(B,\underline{A})$ in the \'etale (or classical)
topology, and $g:X=A^c\lra B$ the torsor constructed from $c$. 

Then the following assertions hold:

(i) Assume that $X$ possesses a quasi-symplectic
structure $\Lambda_X$ such that
$g$ is a Lagrangian fibration with respect to $\Lambda_X$. Then
$A$ possesses also a quasi-symplectic structure $\Lambda_A$ for
which $f$ is Lagrangian, and for any local section $s :V\lra X$ of $g$ inducing
the isomorphism $i_s:X_V=g^{-1}(V)\isoto A_V=f^{-1}(V)$, there exists a
$2$-form $\omega_V\in
H^0(V,\Omega^2)$ such that $\Lambda_X-i_s^*\Lambda_A=g^*\omega_V$
on $V$. Say that such quasi-symplectic structures $\Lambda_X,\Lambda_A$
are associated to each other.

(ii) If $\Lambda_X,\Lambda_A$ are associated quasi-symplectic structures,
then any pair $\Lambda_{1,X}\in
 \Lambda_{X} +g^*H^0(B,\Omega^2)$,
$\Lambda_{1,A}\in
 \Lambda_{A} +f^*H^0(B,\Omega^2)$
is also a pair of associated quasi-symplectic structures.
Moreover, there always exists such a $\Lambda_{1,A}$
which is symplectic.

In particular, if $H^0(B,\Omega^2)=0$ and if there exists a pair
 $\Lambda_X,\Lambda_A$ of associated quasi-symplectic structures,
then it is unique and $\Lambda_A$ is in fact symplectic.

$\Lambda_A$ associated to a given $\Lambda_{X}$ is 
also unique  and symplectique if we impose the
additional condition that the restriction of $\Lambda_A$
to the section of neutral elements of the fibers of $f$ is
identically~$0$.

(iii) Conversely, assume that $A$ possesses a quasi-symplectic structure $\Lambda_A$ such that
$f$ is a Lagrangian fibration with respect to $\Lambda_A$.
Assume also
that at least one of the following two conditions is verified:
$H^1(B,\Omega^2)=0$, or $c$ is of finite order. Then $X$ has a quasi-symplectic structure $\Lambda_{X}$ associated to $\Lambda_A$. In the case when
$c$ is of finite order, $\Lambda_X$ can be chosen symplectic.

(iv) In the algebraic category,
$H^1_{\mbox{\tiny \'et}}(B,\underline{A})$
is torsion, and any torsor with a quasi-symplectic Lagrangian
structure has also a symplectic one.

\end{lemma}

\begin{proof}
For (i), (ii) see Propositions 2.3, 2.6 in \cite{Ma}. In fact, the
assumptions of these propositions are formulated there for symplectic
structures, but the proofs remain valid also for quasi-symplectic ones.
Proposition 2.6 in loc. cit. treats the case of part (iii)
when $H^1(B,\Omega^2)=0$.
We need to introduce some settings of this proof to explain
the case when $c$ is of finite order.

For any $b\in B$ and for a sufficiently small neighborhood $V$
of $b$, there exist local parameters $u_1,\ldots ,u_n$ on $V$
and a local section $s_V$ of $g_V:X_V\lra V$ identifying $X_V$ with $A_V$
in such a way that $s_V$ becomes the section $e$ of neutral elements 
of the fibers of $f$.
One can realize $A_V$ as the quotient $\Omega^1_V/\LLL$ of the
cotangent bundle of $V$ by a
local system $\LLL$ generically of rank $2n$.
Let $z_1,\ldots ,z_n$ be the coordinates on $\Omega^1_{B}$,
dual to $du_1,\ldots ,du_n$. Then the restriction
of  $\Lambda_A$ to $A_V$ can be written
in the form
\begin{equation}\label{eq2}
\Lambda_X=\sum_pdu_p\wedge dz_p+\beta,\;\;\; \beta
=\sum_{p,q}\beta_{pq}(u)du_p\wedge du_q .
\end{equation}
Here $\beta$ does not depend on $z$, as was proved in loc. cit.
Let $\Lambda_{0,A}=\sum_pdu_p\wedge dz_p
=\Lambda_X-\beta$. The invariance of $\Lambda_{0,A}$ under the translations
by sections from $\LLL$ is equivalent to
\begin{equation}\label{eq3}
d(\sum_p\gamma_pdu_p)=0\ ,
\end{equation}
where $\gamma =(\gamma_1,\ldots ,\gamma_n)$ 
is any local section of $\cal L$.
Hence we have also the invariance
under translations by sections from $\QQ\otimes\LLL$, and the finite order
case follows.

The part (iv) is obtained by the successive application of (i), (iii).
\end{proof}

Now we can prove

\begin{theorem}\label{closedness}
There exists a symplectic structure $\tilde{\wedge}_{\mathfrak M}$ on 
${\mathfrak M}$ such that the following two conditions are verified:

(a) $\tilde{\wedge}_{\mathfrak M}-\wedge_{\mathfrak M}\in
H^0(U,p^*\Omega^2_U)$, hence $p$ is Lagrangian w. r. t.
$\tilde{\wedge}_{\mathfrak M}$, and $\tilde{\wedge}_{\mathfrak M}$
also satisfies the conclusion of Theorem \ref{relat}.

(b) $\tilde{\wedge}_{\mathfrak M}$ extends to a symplectic
structure $\Lambda_{\JJB^c}$ associated to $\Lambda_\JJB$.
\end{theorem}

\begin{proof}
We are in the algebraic situation, so we can apply (iii)
to construct a symplectic structure $\Lambda_{\JJB^c}$
associated to the Donagi--Markman form $\Lambda_{\JJB}$.
Define  $\tilde{\wedge}_{\mathfrak M}=\Lambda_{\JJB^c}|_{\mathfrak M}$,
so that (b) is verified.
The definition of the relation ``associated" in (i) of Lemma
\ref{torsors} and Theorem \ref{relat} imply also~(a).
\end{proof}

\begin{remark}
We conjecture that $\wedge_{\mathfrak M}$ is symplectic itself.
\end{remark}

\begin{remark}
The \'etale cohomology class $c$ defining $\JJB^c$ is of order 3,
because one can send the curve $3C$ into  $J^2(X)$,
where $C\subset X=Y\cap H$ is
a normal elliptic quintic,  in a canonical way
in using five times the plane section of $Y$ as the reference curve
for the definition of the Abel--Jacobi map.
So, besides the trivial torsor $\JJB$ and $\JJB^c$, there exists
one more interesting torsor $\JJB^{c^2}$. What geometric interpretation
can one give to this torsor?
\end{remark}

\begin{remark}
One can extend our construction to a bigger open set $\bar{U}$, such that
$\codim \check\PP^5\setminus \bar{U}\geq 2$. We will have 
the same Theorems \ref{relat}, \ref{closedness} for the (quasi)-symplectic
structures $\Lambda_{\bar{\JJB}}$,
${\Lambda}_{\bar{\JJB^c}}$, ${\Lambda}_{\bar{\mathfrak M}}$. As
$H^0(\bar{U},\Omega^2)=0$,
the associated quasi-symplectic structures on the extended families
$\bar{h}:\bar{\JJB}\lra\bar{U}$,
$\bar{h}^c:\bar{\JJB}^c\lra\bar{U}$ are unique and symplectic.
\end{remark}

\begin{proof}
Take for $\bar{U}$ the locus of hyperplanes $H$ such that $H\cap Y$ has at most
one nondegenerate double point as singularity. Any cubic $X$ with only one
isolated singular point contains a normal elliptic quintic curve which does not pass
through the singular point. To see this, take a hyperplane section of $X$
which is a nonsingular cubic surface $S\subset \PP^3$. It contains a $5$-dimensional linear
system of quintics of arithmetic genus 1, and its generic member is a non-singular
space elliptic curve $C$. By Lemma 2.7, c) in \cite{Ma-Ti}, whose proof uses only
the nonsingularity of $X$ at points of $S$, we have $h^0(\NNN_{C/X})=10,h^1(\NNN_{C/X})=0$.
Hence the Hilbert scheme of $X$ is smooth of dimension 10 at $[C]$. The family of quintics
of genus 1 in hypeplane sections of $X$ being 9-dimensional, we see that $C$ is deformable
to a non-space curve, which is the wanted normal elliptic quintic.

Define now $M_X$ to be the subset in the moduli space of sheaves on $X$
consisting of the vector bundles of rank 2 obtained by Serre's construction
from all the normal elliptic quintics in $X$ that do not pass through the
singular point of $X$, and $\bar{\mathfrak M}$ the union of the $M_X$
over all $X=Y\cap H$ with $H\in\bar{U}$.
All our proofs concerning the
dimension, the nonsingularity of $M_X$, ${\mathfrak M}$, as well as
the definition and the proof of the nondegeneracy of $\Lambda_{\mathfrak M}$
work without any modification with $\bar{\mathfrak M}$ in place of
${\mathfrak M}$ for all $H\in\bar{U}$. Thus we obtain the Yoneda
quasi-symplectic structure $\Lambda_{\bar{\mathfrak M}}$
on~$\bar{\mathfrak M}$.

Donagi--Markman \cite[Example 8.22]{D-M}
remark that their symplectic structure extends to a 2-form $\bar{\Lambda}_{\bar{\JJB}}$ on
the generalized relative intermediate Jacobian $\bar{h}:\bar{\JJB}\lra\bar{U}$.
They write that it is not difficult to see that it
is nondegenerate over the boundary. This also follows easily
from our setting.

Indeed, the nondegeneracy of $\Lambda_{\bar{\mathfrak M}}$
implies that the vertical tangent bundle of $\bar{\mathfrak M}$
is isomorphic to $\Omega^1_{\bar{U}}$, hence its associated
group scheme $\bar{\JJB}$ is isomorphic to $\Omega^1_{\bar{U}}/\LLL$,
as in the proof of Lemma \ref{torsors}. The standard
symplectic structure of the cotangent bundle of $\bar{U}$
descends to the quotient by $\LLL$ over $U$, hence also
over $\bar{U}$, as the descent condition
(\ref{eq3}) extends by continuity. Thus we obtain a symplectic structure
$\Lambda_{\bar{\JJB}}$ on $\bar{\JJB}$. As the restrictions
of $\Lambda_{\bar{\JJB}}$, $\bar{\Lambda}_{\bar{\JJB}}$ to $U$
are associated to the same quasi-symplectic structure
on $\JJB^c$,
they differ by a lift of a 2-form on $U$. Hence they define
the same isomorphism between $\Omega^1_{\bar{U}}$ and
$\TTT_{\bar{\JJB}/\bar{U}}$ over $U$, which is the
identity one for $\Lambda_{\bar{\JJB}}$. By continuity,
this holds also over $\bar{U}$. Hence $\bar{\Lambda}_{\bar{\JJB}}$
is nondegenerate. As $H^0(\bar{U},\Omega^2)=0$, we have the uniqueness
over $\bar{U}$, hence
$\bar{\Lambda}_{\bar{\JJB}}=\Lambda_{\bar{\JJB}}$.

The extension of Theorem \ref{relat} follows trivially from
the following observation: if the regular section $\Lambda_{\bar{\mathfrak M}}-
\Psi_{\UUU}^*\Lambda_{\bar{\JJB}}$ of the vector bundle
$\Omega^2_{\bar{\mathfrak M}}$
over $\UUU$
belongs to the vector subbundle
$\bar{p}^*\Omega^2_{\bar{U}}$ over $U\cap\UUU$, then it does over the whole
of $\UUU$.
\end{proof}

\end{document}